\documentclass[12pt]{amsart}
\usepackage{amssymb,amsthm,rotating, mathtools}
\usepackage[all]{xy}

\numberwithin{equation}{section}
 \newtheorem{thm}[equation]{Theorem}

 \newtheorem{defn}[equation]{Definition}
 \newtheorem{prop}[equation]{Proposition}

\newtheorem{cor}[equation]{Corollary}
 \newtheorem{lemma}[equation]{Lemma}

 \theoremstyle{definition}

 \newtheorem{remark}[equation]{Remark}

 \newtheorem{example}[equation]{Example}

\DeclareMathOperator{\Ext}{Ext}
\DeclareMathOperator{\Hom}{Hom}
\DeclareMathOperator{\codim}{codim}
\DeclareMathOperator{\GL}{GL}
\DeclareMathOperator{\Gl}{GL}

\DeclareMathOperator{\Ima}{im}
\DeclareMathOperator{\Proj}{Proj}
\DeclareMathOperator{\sgn}{sgn}
\DeclareMathOperator{\Span}{Span}
\DeclareMathOperator{\Sym}{Sym} 
\DeclareMathOperator{\vol}{vol}
\DeclareMathOperator{\diag}{diag}

\DeclareMathOperator{\U}{U} 

\newcommand{\sta}{\diamond}

\newcommand{\DOT}{\setlength{\unitlength}{1pt}\begin{picture}(2.5,2)
                  (1,1)\put(2,3.5){\circle*{2}}\end{picture}}
\newcommand{\coh}{{\rm H}}
\newcommand{\del}{\partial}
\newcommand{\HH}{{\rm HH}}
\newcommand{\HHD}{{\rm HH}^{\DOT}}
\newcommand{\Wedge}{\textstyle\bigwedge}
\newcommand{\CC}{\mathbb{C}}
\newcommand{\NN}{\mathbb{N}}
\newcommand{\ZZ}{\mathbb{Z}}
\newcommand{\ot}{\otimes}

\newcommand{\ep}{\epsilon}

\newcommand{\scooch}{\hspace{-.25ex}}
\newcommand{\ignore}[1]{\relax}

\newcommand{\ta}{\Upsilon}

\title[Finite groups acting linearly]
{Finite groups acting linearly: Hochschild cohomology and the cup product}

\date{November 4, 2009}
\author{Anne V.\ Shepler}
\address{Department of Mathematics, University of North Texas,
Denton, Texas 76203, USA}
\email{ashepler@unt.edu}
\author{Sarah Witherspoon}
\address{Department of Mathematics\\Texas A\&M University\\
College Station, Texas 77843, USA}\email{sjw@math.tamu.edu}
\thanks{The first author was partially supported by NSF grants
\#DMS-0402819 and \#DMS-0800951 and a research fellowship from the 
Alexander von Humboldt Foundation. 
The second author was partially supported by NSA grant
\#H98230-07-1-0038 and NSF grant
\#DMS-0800832.
Both authors were jointly supported by Advanced Research Program Grant 010366-0046-2007
from the Texas Higher Education Coordinating Board.}

\begin{document}
\maketitle

\begin{abstract}
When a finite group acts linearly on a complex vector space, the 
natural semi-direct product of the group and the polynomial ring over the
space forms a skew group algebra.  
This algebra plays the role of the coordinate ring of the
resulting orbifold
and serves as a substitute for the ring of invariant polynomials
from the viewpoint of geometry and physics.
Its Hochschild cohomology 
predicts various Hecke algebras and deformations of the orbifold.
In this article,
we investigate the ring structure of the Hochschild cohomology of the skew
group algebra.
We show that the cup product coincides with a natural smash product,
transferring the cohomology of a group action
into a group action on cohomology.  
We express the algebraic structure of Hochschild cohomology in terms of a partial order
on the group (modulo the kernel of the action).  This partial order arises after
assigning to each group element the codimension
of its fixed point space. 
We describe the algebraic structure 
for Coxeter groups, where this partial order is given by
the reflection length function; a similar combinatorial description 
holds for an infinite family of complex reflection groups.
\end{abstract}

\section{Introduction}

Physicists often regard space as a Calabi-Yau manifold $M$ endowed
with symmetries forming a group $G$.  
The orbifold $M/G$ is not smooth in general,
and they regularly shift focus from the orbifold 
to a desingularisation and its coordinate ring.
In the affine case, we take $M$ to be a finite dimensional, complex vector
space $V$ upon which a finite group $G$
acts linearly.
The orbifold $V/G$ may then be realized as an algebraic variety
whose coordinate ring is the ring of invariant
polynomials $S(V^*)^G$ on the dual space $V^*$ (see Harris~\cite{Harris}). 
The variety $V/G$ is nonsingular exactly when the action of $G$ on $V$ 
is generated by reflections.
When the orbifold $V/G$ is singular,
we seek to replace
the space of invariant polynomials with a natural
algebra attached to $V/G$ playing the role of a coordinate ring.

In situations arising naturally in physics,
a noncommutative substitute for the commutative invariant ring $S(V^*)^G$
is provided by
the skew group algebra $$S(V^*)\# G=S(V^*)\rtimes G,$$
the natural semi-direct product of $G$ with the symmetric algebra
$S(V^*)$.
One resolves the singularities of $V/G$ with a smooth Calabi-Yau variety
$X$ whose coordinate ring behaves
as the skew group algebra.
Indeed, the McKay equivalence implies (in certain settings) that Hochschild cohomology
sees no difference:
$$
\HHD(S(V^*)\#G) \cong \HHD(\mathcal{O}_X)
$$
as algebras under both the cup product and
Lie bracket.
Note that the Hochschild cohomology $\HHD(S(V^*)\#G)$ can also be used 
to recover the orbifold (or stringy) cohomology $\coh_{orb}^{\DOT}(V/G)$
(which is isomorphic to the singular cohomology of $X$).
(See~\cite{BezrukavnikovKaledin}, \cite{BridgelandKingReid}, 
\cite{CaldararuGiaquintoWitherspoon}, 
\cite{GinzburgKaledin}, \cite{Keller}.)
The cohomology and deformations of $S(V)\# G$
appear in various other areas of mathematics as well---for example,
combinatorics, representation theory,
Lie theory, noncommutative algebra, and invariant theory (see, for example,
Etingof and Ginzburg~\cite{EtingofGinzburg}).

In this paper, we consider any finite group $G$ acting linearly on $V$
and explore the rich algebraic structure of the Hochschild cohomology
of $S(V)\# G$ under the cup product.
This structure is interesting
not only in its own right, but 
also because of possible applications in algebra and representation theory.
For example, the graded Lie bracket, which predicts
potential deformations (like sympletic reflection algebras and graded Hecke
algebras), 
is a graded derivation on Hochschild
cohomology with respect to the cup product. 
The representation theory of finite dimensional algebras
provides an application of the cup product in a different setting:
Often, one may associate an algebraic variety to each module over the algebra
using the ring structure of its Hochschild cohomology;
the collection of such varieties 
provide a coarse invariant of the representation theory of the algebra
(see, for example, Snashall and Solberg~\cite{SnashallSolberg}).

For any algebra $A$ over a field $k$, Hochschild cohomology $\HH^{\DOT}(A)$
is the space $\Ext^{\DOT}_{A\ot A^{op}}(A,A)$, which
is a Gerstenhaber algebra under the two compatible operations, cup product and
bracket. 
Both operations
are defined initially on the bar resolution, a natural $A\ot A^{op}$-free
resolution of $A$. 
For $A=S(V)\# G$, an explicit description of Hochschild cohomology 
$\HHD(A)$ arises not from the bar resolution of $A$, but instead 
from a Koszul resolution:
$ \HH^{\DOT}(S(V)\# G)$ is isomorphic
to the $G$-invariant subalgebra of $\HH^{\DOT}(S(V), S(V)\# G)$,
which is computed with a Koszul resolution of $S(V)$. 
We analyze the cup product
on $\HHD(S(V)\# G)$ by taking advantage of these two different manifestations of
cohomology arising from two different resolutions.
We show that the cup product on $\HH^{\DOT}(S(V)\# G)$ may be written as the
cup product on $S(V)$ twisted by the action of the group.
This perspective yields convenient descriptions of the
ring structure of cohomology.
(We study the Gerstenhaber bracket in a future work.)

Over the real numbers,
the cup product on $\HH^{\DOT}(S(V)\# G)$ has been studied 
in related settings (e.g., see~\cite{PPTT}). 
But note that this analysis of Hochschild cohomology 
does not readily extend to our setting.
In this paper, we do not assume that $G$ acts symplectically or even
faithfully,  as there are interesting
applications in which it does not
(see~\cite{Chmutova,Dezelee,SheplerWitherspoon}).
Over the real numbers, $V$ and $V^*$ are naturally $G$-isomorphic, 
which may simplify several aspects of the theory.
We develop the theory in the richer setting of complex affine space.
Note that our results actually hold over any field containing the
eigenvalues of the action of $G$ on $V$ in which $|G|$ is invertible.
Anno~\cite{Anno} also gave a cup product formula in the geometric setting
over fairly general fields;
we give a natural interpretation of the resulting ring structure
from a purely algebraic and combinatorial point of view.

In Section~\ref{lemmas}, we define the ``codimension poset'' which arises from
assigning to each group element the codimension of its fixed point space.
We also posit a few observations needed later on the geometry of finite group
actions.
In Sections~\ref{cohomologydefs} and~\ref{koszul}, we establish definitions
and notation and recall necessary facts about the bar and Koszul resolutions. 
We review the structure
of $\HHD(S(V)\# G)$ as a graded vector space as well, originally found independently by
Farinati~\cite{Farinati} and Ginzburg and Kaledin~\cite{GinzburgKaledin} for
faithful group actions.
In Section~\ref{Sec:tau}, we define
a combinatorial map in terms of Demazure operators and quantum differentiation
which allows for conversion between complexes.  This
combinatorial conversion map (introduced in~\cite{SWchain})
induces isomorphisms on cohomology.

In Sections~\ref{smashstructure} and~\ref{cup=smash}, we transform the
``cohomology of a group action'' into a ``group action on cohomology''.
One may first take Hochschild cohomology and
then form the skew group algebra, or one may reverse this order of operation.
In Section~\ref{smashstructure}, we compare resulting algebras for $S(V)$:
$$
\HH^{\DOT}(S(V)\# G)
\quad\text{versus}\quad 
\HH^{\DOT}(S(V))\# G
\quad\text{versus}\quad 
\HHD(S(V))\ot \CC[G]\ .
$$
Since $\HH^{\DOT}(S(V)\# G)$ is 
the $G$-invariant subalgebra of $\HH^{\DOT}(S(V), S(V)\# G)$,
we focus on this latter ring.
We show that the smash product on $\HH^{\DOT}(S(V))\# G$
induces a smash product on $\HH^{\DOT}(S(V), S(V)\# G)$,
which we then 
view as an algebra under three operations:
\begin{itemize}
\item the cup product $\smile$ induced from the bar resolution of $S(V)$, 
\item the smash product $\sta$ induced from  $\HHD(S(V))\# G$, and
\item the usual multiplication in the tensor algebra product $\HHD(S(V))\ot
  \CC[G]$.
\end{itemize}

In Section~\ref{cup=smash} (see Theorem~\ref{thm:cup=smash}), 
we show that these three algebraic operations {\em coincide}.
This yields a simple formula 
for the cup product on $\HHD(S(V),S(V)\# G)$ (see Theorem~\ref{cupformula};
cf.\ Anno~\cite{Anno})
and implies that $\HHD(S(V)\# G)$ is isomorphic to an {\em algebra
subquotient} of $\HHD(S(V))\# G$ (see Corollary~\ref{algebrasubquot}).  
These results express the 
cup product on $\HHD(S(V)\#G)$ (at the cochain level) 
as the cup product on $\HHD(S(V))$ twisted by the group $G$.

In Section~\ref{volumealgebra}, we identify an interesting 
graded subalgebra (the ``volume subalgebra'') of $\HHD(S(V),S(V)\# G)$
whose dimension is the order of $G$. When
$G$ is a subgroup of the symplectic group $\text{Sp}(V)$, its $G$-invariant
subalgebra is isomorphic to the cohomology of the $G$-invariant
subalgebra of the Weyl algebra.
In this case,
it is also isomorphic to the orbifold cohomology of $V/G$.
(See Remark~\ref{describesubalg}.)
We thus display the orbifold cohomology $\coh_{\text{orb}}^{\DOT}(V/G)$ 
as a natural {\em subalgebra} of $\HH^{\DOT}(S(V)\#G)$.

In Section~\ref{poset}, we describe generators of cohomology (as an algebra)
via the codimension poset.
The partial order is defined on $G$ modulo the kernel of its action on $V$.
We view division in the volume algebra
as a purely geometric construction by interpreting results in terms 
of this poset.
Generators for the Hochschild cohomology $\HHD(S(V),S(V)\#G)$
arise from minimal elements in the poset (with the identity removed);
see Corollaries~\ref{generatednonfaithful} and~\ref{generatedfaithful}.

In Section~\ref{reflectiongroups}, we explore implications for
reflection groups.  The reflection length of a group element is the length of
a shortest word expressing that element as a product of reflections.
For Coxeter groups and many complex reflection groups,
reflection length coincides with the codimension of the fixed point space.
In this case, the codimension poset
appears as a well-studied poset (arising from reflection length) in the
theory of reflection groups.
We show (in Corollary~\ref{generatorsforCoxeter}) 
that for Coxeter groups and many other complex
reflection groups, the Hochschild cohomology $\HHD(S(V),S(V)\# G)$ 
is generated as an algebra in homological degrees $0$ and $1$,
in analogy with the Hochschild-Kostant-Rosenberg Theorem
for smooth commutative algebras (such as $S(V)$).

Finally, in Section~\ref{last-section}, we return to $\HHD(S(V )\# G)$ 
as the $G$-invariant subalgebra of $\HHD(S(V),S(V)\# G)$.
We point out the simple cup product structure in cohomological degrees 0 and
1.  In Theorem~\ref{prop:formula},
we use standard group-theoretic techniques (Green/Mackey functors and transfer
maps) to describe the product
on $\HHD(S(V)\# G)$.  We use the natural  symplectic action 
of the symmetric group to give a nontrivial example.

\vspace{3ex}

\section{Poset and volume forms}\label{lemmas}

We begin by collecting several geometric observations needed later.
Let $G$ be a finite group and $V$ a $\CC G$-module 
of finite dimension $n$.
Denote the image of $v\in V$ under the action of 
$g\in G$ by $ {}^gv$.  We work with the induced group action on maps:
For any function $\theta$ and any element $h\in \Gl(V)$ acting on
its domain and range, 
we define the map ${}^h\theta $
by $(^h\theta)(v) := \hphantom{\,}^h(\theta(^{h^{-1}}v))$.
Let $V^*$ denote the contragredient (or dual) representation.
For any basis $v_1,\ldots,v_n$ of $V$, let
$v_1^*,\ldots,v_n^*$ be the dual basis of $V^*$.
Let $V^G$ be the set of $G$-invariants in $V$:
$V^G=\{v \in V: {}^{g}v=v \text{ for all } g\in G\}$. 
For any $g\in G$, 
set $V^g = \{v\in V: {}^gv=v\}$, the fixed point set of $g$ in $V$.
Since $G$ is finite, we assume $G$ acts by isometries on $V$ (i.e.,
$G$ preserves a given inner product).
 All tensor and exterior products will be taken
over $\CC$ unless otherwise indicated.

We regard $G$ modulo the kernel of its action on $V$ 
as a poset using the following lemma, which  is no surprise
(see~\cite{BradyWatt} or~\cite{GinzburgKaledin}, for example). 
Note that $(V^g)^{\perp} = \Ima (1-g)$ for all $g$ in $G$.
\vspace{2ex}
\begin{lemma}\label{codims}
Let $g,h\in G$. The following are equivalent:

(i) $(V^g)^{\perp}\cap (V^h)^{\perp} = 0$,

(ii) $V^g+V^h = V$,

(iii) $\codim V^g + \codim V^h = \codim V^{gh}$,

(iv) $(V^g)^\perp\oplus(V^h)^\perp=(V^{gh})^\perp$.

\noindent
Any of these properties implies that 
 $V^g\cap V^h = V^{gh}$.
\end{lemma}
\vspace{1ex}
\begin{proof}
Taking orthogonal complements yields
the equivalence of (i) and (ii). 

Assume
(ii) holds and write $u\in V^{gh}$
as $v+w$ where $v\in V^g$ and $w\in V^h$.
Then ${}^hv+{}^hw={}^{g^{-1}}v+{}^{g^{-1}}w$ and
${}^hv+w=v+{}^{g^{-1}}w$, i.e., ${}^{(h-1)}v={}^{(g^{-1}-1)}w$.
But $(V^h)^{\perp}=\Ima(h-1)$ while $(V^g)^{\perp}=(V^{g^{-1}})^{\perp}
=\Ima(g^{-1}-1)$, and the intersection of these spaces is 0 by (i),
hence ${}^{(h-1)}v=0={}^{(g^{-1}-1)}w$. Therefore, $v\in V^h$, $w\in V^g$,
and $u \in V^g \cap V^h$, and thus
$V^{gh}\subset V^g\cap V^h$.  The reverse inclusion 
is immediate, hence $V^{gh}=V^g\cap V^h$.
We take orthogonal complements and observe that $(V^{gh})^{\perp} = 
(V^g)^{\perp}\oplus (V^{h})^{\perp}$. 
A dimension count then gives (iii).

To show (iii) implies (iv), note
$(V^{gh})^{\perp}\subset (V^g)^{\perp} + (V^h)^{\perp}$
since $V^g\cap V^h\subset V^{gh}$.  
By (iii), this containment is forced to be an equality and
the sum is direct: $(V^g)^\perp\oplus(V^h)^\perp=(V^{gh})^\perp$.
As (iv) trivially implies (i), we are finished.
\end{proof}
\vspace{2ex}


Let $K$ be the kernel
of the representation of $G$ acting on $V$: $$K:=\{k\in G: V^k=V\}.$$

\vspace{3ex}
\begin{defn}\label{binaryrelation}
Define a binary relation $\leq$ on $G$ 
by $g \leq h$ whenever
$$
  \codim(V^g) + \codim (V^{g^{-1}h}) = \codim (V^h).
$$
By Lemma~\ref{codims}, this codimension condition holds exactly when
$$(V^g)^\perp\oplus(V^{g^{-1}h})^\perp=(V^{h})^\perp\ .$$
This induces a binary relation $\leq$ on the quotient group $G/ K$ 
as well: For $g,h$ in $G$, define $gK\leq hK$ when $g \leq h$. Note the relation does
not depend on choice of representatives of cosets, 
as $V^{g}=V^{h}$ whenever $gK=hK$ for $g,h$ in $G$.
\end{defn}
\vspace{3ex}

The relation $\leq$ appears in work of Brady and Watt~\cite{BradyWatt}
on orthogonal transformations. Their arguments apply equally well 
to our setting of isometries with respect to some
inner product and the quotient group $G/K$.  (Note that if $G$ does not act 
faithfully, the binary relation $\leq$ on $G$ may not be anti-symmetric
and thus may not define a partial order on $G$.)
\vspace{2ex} 
\begin{lemma}[Brady and Watt~\cite{BradyWatt}]\label{posetlemma}
The relation $\leq$ is a partial order on $G/K$.
If $G$ acts faithfully, the relation $\leq$ is a partial order on $G$.
\end{lemma}

\vspace{3ex}

We shall use the following elements in the sequel.

\begin{defn}
For each $g\in G$, let $\vol^{\perp}_g$ be a choice of nonzero
element in the one-dimensional space $\Wedge^{\codim V^g}((V^g)^{\perp})^*$.
\end{defn}
\vspace{3ex}

We show in the next lemma how these choices determine a multiplicative
cocycle.
A function $\vartheta: G\times G \rightarrow \CC$ is a 
{\bf multiplicative $2$-cocycle} on $G$
if 
$$
   \vartheta(gh,k)\vartheta(g,h) = \vartheta(g,hk)\vartheta(h,k)
$$
for all $g,h,k\in G$.
We may use any such cocycle 
to define a new algebra structure on the
group algebra $\CC[G]$, a generalization of a twisted group algebra
(in which the values of $\vartheta$ may include 0):
Let $\CC^{\vartheta}[G]$ be the $\CC$-algebra with basis $G$ and
multiplication $g\cdot_{\vartheta} h = \vartheta(g,h)\, gh$ for all $g,h\in G$.
Associativity is equivalent to the 2-cocycle identity.
If $\vartheta(g,1_G)=1=\vartheta(1_G,g)$ for all $g\in G$, where $1_G$ denotes
the identity element of $G$, then $\CC^{\vartheta}[G]$
has multiplicative identity $1_G$.

We canonically embed each space $\Wedge ((V^g)^{\perp})^*$ into $\Wedge V^*$.

\vspace{2ex}
\begin{prop}\label{wedgeofvols}
For all $g$ and $h$ in $G$,
$$
\vol_g^\perp\wedge \vol_h^\perp=\vartheta(g,h)\vol_{gh}^\perp
$$ 
in $\Wedge V^*$ where 
$\vartheta: G\times G\rightarrow \CC$ is a (multiplicative) 2-cocycle on $G$
with
$$\vartheta(g,h)\neq 0
\quad\text{if and only if}\quad
g\leq gh \ .
$$
Under wedge product, the algebra $\Wedge\{\vol_g^\perp: g\in G\}$ 
is isomorphic to the (generalized) twisted group 
algebra $\CC^{\vartheta}[G]$.
\end{prop}
\vspace{1ex}
\begin{proof}
Let $g,h$ be any pair of elements in $G$.
Then $V^g \cap V^h \subset V^{gh}$, and hence 
$(V^{gh})^\perp \subset (V^g)^\perp + (V^h)^\perp$.
If the sum is direct, then by Lemma~\ref{codims},
we have equality of vector subspaces:
$(V^{gh})^\perp = (V^g)^\perp \oplus (V^h)^\perp$.
If the sum is not direct, then 
$\vol_g^\perp\wedge\vol_h^\perp=0$.
In either case,  the product
$\vol_g^\perp\wedge \vol_h^\perp$ is a (possibly zero)
scalar multiple of $\vol_{gh}^\perp$. 
Hence, there is a scalar $\vartheta(g,h)\in \CC$ such that
$$
   \vol^{\perp}_g\wedge \vol^{\perp}_h = \vartheta(g,h)\vol^{\perp}_{gh}.
$$
Note that $\vartheta(g,h)\neq 0$ if and only if $\codim V^g + \codim V^{h} =
   \codim V^{gh}$
(by Lemma~\ref{codims}) if and only if $g\leq gh$.
By associativity of the exterior algebra, the 
function $\vartheta: G\times G\rightarrow \CC$ is a (multiplicative) 2-cocycle on $G$. 
\end{proof}


We shall also need the following easy lemma, which is a consequence of 
the fact that $(V^{g})^\perp=\Ima (1-g)$ for all $g$:
\vspace{2ex}
\begin{lemma}\label{wedgewithperpspace}
For all $g$ in $G$:
\begin{itemize}
\item
$f -\, ^gf$ lies in the ideal $I((V^g)^{\perp})$ for all $f \in S(V)$, and
\item
$dv\wedge \vol_g^{\perp}=\, ^gdv\wedge \vol_g^{\perp}$
\ \ for all $dv \in \bigwedge^{\DOT} V$.
\end{itemize}
\end{lemma}
\vspace{3ex}

\section{Skew group algebra and Hochschild cohomology}
\label{cohomologydefs}

In this section, we recall the basic definitions 
of the skew group algebra and Hochschild
cohomology, as well as a fundamental theorem describing the cohomology
of the skew group algebra as a space of invariants.
We work over the complex numbers $\CC$.

Let $A$ denote any $\CC$-algebra 
on which $G$ acts by automorphisms.
The {\bf skew group algebra} (or {\bf smash product}) $A\# G$ is the vector space
$A\ot \CC G$ with multiplication
given by
$$
  (a\ot g)(b\ot h)=a ( {}^{g}b) \ot gh
$$
for all $a,b\in A$ and $g,h\in G$.

The {\bf Hochschild cohomology} of a $\CC$-algebra $A$ (such as $S(V)$ or $S(V)\# G$),
with coefficients in an $A$-bimodule $M$,
is the graded vector space
$\HH^{\DOT}(A,M)=\Ext^{\DOT}_{A^e}(A,M)$, where $A^e=A\ot A^{op}$ acts
on $A$ by left and right multiplication.
We abbreviate $\HH^{\DOT}(A)=\HH^{\DOT}(A,A)$.

To construct the cohomology $\HHD(A,M)$, one applies the
functor $\Hom_{A^e}( - ,M)$ to a projective resolution of $A$ as an $A^e$-module,
for example, to the {\bf bar resolution}
\begin{equation}\label{barcomplex}
  \cdots\stackrel{\delta_3}{\longrightarrow} A^{\ot 4}
\stackrel{\delta_2}{\longrightarrow} A^{\ot 3}
\stackrel{\delta_1}{\longrightarrow} A^e
\stackrel{m}{\longrightarrow} A \rightarrow 0,
\end{equation}
where $\delta_p(a_0\ot\cdots\ot a_{p+1}) = \sum_{j=0}^p (-1)^j a_0
\ot\cdots\ot a_ja_{j+1}\ot\cdots\ot a_{p+1}$ for all $a_0,\ldots, a_{p+1}\in A$,
and $m$ is
multiplication.
For all $p\geq 0$, $\Hom_{A^e}(A^{(p+2)}, M)\cong \Hom_{\CC}(A^{\ot p},M)$,
and we identify these two vector spaces in what follows.

When $M$ is itself an algebra, 
the Hochschild cohomology $\HH^{\DOT}(A,M)$ is a graded associative
algebra under the {\bf cup product}, 
defined at the cochain level on the bar complex (see~\cite[\S7]{Gerstenhaber63}): 
Let $f\in\Hom_{\CC}(A^{\ot p},M)$ and $f'\in\Hom_{\CC}(A^{\ot q},M)$;
then the cup product
 $f\smile f'$ in $\Hom_{\CC}(A^{\ot (p+q)},M)$
is given by
$$
  (f\smile f')(a_1\ot\cdots\ot a_{p+q}) = f(a_1\ot\cdots\ot a_p)
       f'(a_{p+1}\ot\cdots\ot a_{p+q})
$$
for all $a_1,\ldots,a_{p+q}\in A$.
We seek to describe the algebra structure under cup product explicitly in
the case $A=S(V)$, $M=S(V)\# G$, and in the case $A=M=S(V)\# G$.
These two cases are related by a well-known theorem that we state next.

Since $|G|$ is invertible, a result of
\c{S}tefan~\cite[Cor.\ 3.4]{Stefan} implies in our setting that there is a $G$-action 
giving an isomorphism of graded algebras (under the cup product):
\vspace{2ex}
\begin{thm}\label{isos}
$$
\begin{aligned}
  \HH^{\DOT}(S(V)\# G) \ \ 
&\cong \ \ 
\HH^{\DOT}(S(V), S(V)\#G)^G
\ .
\end{aligned}
$$
\end{thm}
\vspace{2ex}
\noindent
(Specifically, the action of $G$ on $V$ extends naturally to the bar complex of $S(V)$ 
and thus induces an action on $\HH^{\DOT}(S(V),S(V)\#G)$, from which we define
$G$-invariant cohomology in the theorem.
In fact, any projective resolution of $S(V)$ compatible with the action of $G$
may be used to define the $G$-invariant cohomology.)
We thus concentrate on describing the cup product on $\HH^{\DOT}(S(V),S(V)\#G)$.
\vspace{3ex}
\section{Koszul and bar resolutions}\label{koszul}

One may use either the Koszul or the bar resolution of $S(V)$ to describe the cohomology 
$\HH^{\DOT}(S(V), S(V)\#G)$
and thus its $G$-invariant subalgebra, $\HH^{\DOT}(S(V)\# G)$. 
The {\bf Koszul resolution} is the following free $S(V)^e$-resolution of $S(V)$:
\begin{equation}\label{koszul-res2}
 \cdots \stackrel{d_3}{\longrightarrow} S(V)^e\ot \Wedge^2V
   \stackrel{d_2}{\longrightarrow} S(V)^e\ot \Wedge^1 V 
  \stackrel{d_1}{\longrightarrow} S(V)^e \stackrel{m}{\longrightarrow} S(V)
   \rightarrow 0, 
\end{equation}
where the  differential $d$ is given by 
\begin{equation}\label{koszul-diff}
d_p(1\ot 1\ot v_{j_1}\wedge\cdots\wedge v_{j_p}) = 
   \sum_{i=1}^p (-1)^{i+1} (v_{j_i}\ot 1 - 1\ot v_{j_i})\ot
   (v_{j_1}\wedge\cdots\wedge \hat{v}_{j_i}\wedge\cdots\wedge v_{j_p}),
\end{equation}
for all $v_{j_1},\ldots, v_{j_p}\in V$.
Let $\Phi$ be the canonical inclusion (a chain map)
of the Koszul resolution~(\ref{koszul-res2})
into the bar resolution~(\ref{barcomplex}):
$$
\xymatrix{
\cdots \ar[r] & S(V)^{\ot 4}\ar[r]^{\delta_2}
               & S(V)^{\ot 3}\ar[r]^{\delta_1}
               & S(V)^e  \ar[r]^m  
               & S(V)  \ar[r] 
               & 0\\
\cdots \ar[r] & S(V)^e\ot \Wedge ^2 V \ar[r]^{d_2}\ar@<-2pt>[u]_{\Phi_2}
               & S(V)^e\ot \Wedge ^1 V \ar[r]^{\hspace{.6cm}d_1}\ar@<-2pt>[u]_{\Phi_1}
               & S(V)^e \ar[r]^{m} \ar@<-2pt>[u]_{=}
               & S(V)\ar[r] \ar@<-2pt>[u]_{=} & 0 ,
}
$$
that is, for all $p\geq 1$,  
\begin{eqnarray*}
   \Phi_p &:& S(V)^e\ot \Wedge^p V\rightarrow S(V)^{\ot (p+2)},
\end{eqnarray*}
\begin{equation}\label{phik}
  \Phi_p(1\ot 1\ot v_{j_1}\wedge\cdots\wedge v_{j_p})
   = \sum_{\pi\in \Sym_p}\sgn(\pi)\ot v_{j_{\pi(1)}}\ot\cdots\ot
    v_{j_{\pi(p)}}\ot 1\ ,
\end{equation}
for all $v_{j_1},\ldots,v_{j_p}\in V$, where $\Sym_p$ denotes the 
symmetric group on the set $\{1,\ldots,p\}$.
Note that $\Phi$ is invariant
under the action of $\GL(V)$, i.e., $^h\Phi = \Phi$ for all $h$ in $\Gl(V)$.

One finds cohomology $\HH^{\DOT}(S(V), M)$ by applying the functor 
$\Hom_{S(V)^e}(-, M)$ to either of the above two resolutions
and dropping the term $\Hom_{S(V)^e}(S(V),M)$.
We make the customary identifications: For each $p$, set
$\Hom_{S(V)^e}(S(V)^{\ot (p+2)},M)=\Hom_{\CC}(S(V)^{\ot p}, M)$ as before, 
and 
\begin{equation}\label{identify}
\Hom_{S(V)^e}(S(V)^e \ot \Wedge^{p} V ,M)
=\Hom_{\CC}(\Wedge^{p} V, M)
=\Wedge^pV^*\ot M\ .
\end{equation}
In the case $M=S(V)\# G$, we
write $S(V)\ot\Wedge^p V^*\ot\CC[G]$
for the vector space $\Wedge^p V^* \ot M$.

We obtain a commutative diagram giving two different
cochain complexes describing the  
cohomology $\HH^{\DOT}(S(V), S(V)\# G)$.  
In~\cite{SWchain}, we introduced 
a ``combinatorial converter'' map $\ta$
(whose definition is recalled  in the next section), which serves as an inverse to the induced  map
$\Phi^*$ and converts between complexes:
\begin{equation}\label{big-diagram-tau}
\xymatrix{
 \Hom_{\CC}(S(V)^{\ot p}, S(V)\# G) \ar[r]^{\delta^*}\ar@<2pt>[d]^{\Phi^*_p}
     & \Hom_{\CC}(S(V)^{\ot (p+1)}, S(V)\# G)\ar@<2pt>[d]^{\Phi^*_{p+1}}\\
  S(V)\ot \Wedge^p V^*\ot \CC[G] \ar[r]^{d^*}  \ar@<2pt>[u]^{\ta_p}
  & S(V)\ot\Wedge^{p+1}V^*\ot\CC[G]\ar@<2pt>[u]^{\ta_{p+1}}
\ \ \ .
}
\end{equation}
We use the maps $\ta$ and $\Phi^*$ in our analysis of the cup product
in later sections.

We describe cohomology explicitly in terms of cocycles and coboundaries.
Under the identification (\ref{identify}), 
Hochschild cohomology $\HHD(S(V), S(V)\#G)$
arises from the complex of cochains
\begin{equation}
\label{cochains}
 C^{\DOT}:=
\bigoplus_{g\in G} S(V)\otimes \Wedge^{\DOT} V^*\ot g \ .
\end{equation}
One may determine the set of cocycles and coboundaries explicitly 
as the kernel and image of the induced map $d^*$ (under the
above identifications (\ref{identify})).
We set
\begin{equation}\label{cocycles}
Z^{\DOT}:=
\bigoplus_{g\in G} S(V) \otimes \Wedge^{\DOT-\codim V^g}
(V^g)^*\otimes\Wedge^{\codim V^g}((V^g)^\perp)^*\ot g\ , 
\end{equation}
a subspace of the space of cocycles, and
$$
B^{\DOT}:=
\bigoplus_{g\in G} I((V^g)^\perp)\otimes \Wedge^{\DOT-\codim V^g}
  (V^g)^*\otimes\Wedge^{\codim V^g}((V^g)^\perp)^* \ot g\ ,
$$
a subspace of the space of coboundaries,
where $I((V^g)^\perp)$ is the ideal of $S(V)$ generated by $(V^g)^\perp$.
(We agree that a negative exterior power of a space is defined to be 0.)
We regard these subspaces $Z^{\DOT}$ and $B^{\DOT}$
as subsets of the cochains $C^{\DOT}$ after making
canonical identifications (we identify $W_1\ot W_2$ with $W_1\wedge W_2$ for any subspaces $W_1,
W_2$ of $V$ intersecting trivially).  We refer to cochains, cocycles, and
coboundaries as vector forms ``tagged'' by the group elements indexing
the direct summands above.

The next remark explains that we may view $Z^{\DOT}$ as a substitute for
the set of cocycles and $B^{\DOT}$ as a substitute for the set of
coboundaries.  We use this alternate description of cohomology 
(as $Z^{\DOT}/B^{\DOT}$) 
in Section~\ref{smashstructure} to expose a smash product structure.


\vspace{3ex}
\begin{remark}\label{H_g}
Note that one may use the isomorphism $V/(V^g)^\perp\cong V^g$
to select a set of a representatives of cohomology classes:
Let
$$
 H^{\DOT}:=\bigoplus_{g\in G}\ S(V^g)\ot  \Wedge^{\DOT - \codim V^g}(V^g)^*
     \ot \Wedge^{\codim V^g}((V^g)^{\perp})^*\ot g\ .
$$
Then $H^{\DOT}\cong Z^{\DOT}/ B^{\DOT}$
via the canonical map
$\Proj_H:Z^{\DOT}/B^{\DOT}\rightarrow H^{\DOT}$ 
induced from the compositions, for each $g$ in $G$,
$$
S(V)\longrightarrow S(V)/ I((V^g)^\perp) \stackrel{\sim}{\longrightarrow} S(V^g)\ .
$$
Farinati~\cite{Farinati}
and Ginzburg and Kaledin~\cite{GinzburgKaledin} showed
that
$$
\HH^{\DOT}(S(V),S(V)\#  G) \cong H^{\DOT}\ 
\cong Z^{\DOT}/ B^{\DOT}.
$$
\end{remark}
\vspace{3ex}

\section{Combinatorial converter map}\label{Sec:tau}

We recall the definition
of the combinatorial converter map $\ta$
(in Diagram~\ref{big-diagram-tau}) introduced in~\cite{SWchain}.
A nonidentity element of $\GL(V)$ is a {\bf reflection}
if it fixes a hyperplane in $V$ pointwise.
Given any basis $v_1, \ldots, v_n$ of $V$, let $\partial /\partial v_i$ denote the usual
partial differential operator with respect to $v_i$.
In addition, given a complex number $\epsilon\neq 1$, we 
define the $\epsilon$-{\bf quantum partial
differential operator} with respect to
$v:=v_i$ as the scaled Demazure
(BBG) operator 
$\del_{v, \epsilon}:S(V)\rightarrow S(V)$
given by
\begin{equation}\label{qpd}
\del_{v, \epsilon}(f)\ = \ (1-\epsilon)^{-1}\ \frac{f-\, ^s\hspace{-.5ex}f}{v}
\ = \ \frac{f-\, ^s\hspace{-.5ex}f}{v-\, ^sv}\ ,
\end{equation}
where $s\in\Gl(V)$ is the reflection whose matrix with respect to the 
basis $v_1, \ldots, v_n$ is 
$\diag(1,\ldots, 1, \epsilon, 1, \ldots, 1)$
with $\epsilon$ in the $i$-th slot.
 Set $\del_{v,\ep}=\del /\del v$ when $\ep=1$.
The operator $\del_{v,\epsilon}$ coincides with the usual definition of
quantum partial differentiation:
One takes the ordinary partial derivative with respect to $v$ but instead
of multiplying each monomial by its degree $k$ in $v$, one multiplies by
the quantum integer
$$
   [k]_{\epsilon} := 1 + \epsilon + \epsilon^2 +\cdots + \epsilon^{k-1}.
$$

We define the map $\ta$ in terms of these Demazure operators.  
For each $g$ in $G$,
fix a basis $B_g=\{v_1,\ldots, v_n\}$
of $V$ consisting of eigenvectors of $g$
with corresponding eigenvalues $\epsilon_1, \ldots, \epsilon_n$.
Decompose $g$ into (commuting) reflections 
diagonal in this basis: 
Let $g=s_1\cdots s_n$ where each $s_i$ is either the identity or a reflection 
defined by
$s_i(v_j) = v_j$ for $j\neq i$ and $s_i(v_i)=\epsilon_i v_i$.
Let $\del_{i}:=\del_{v_i, \epsilon_i}$, 
the quantum partial derivative
with respect to $B_g$.
Recall that $C^{\DOT}$ denotes cochains (see (\ref{cochains})).

\vspace{3ex}
\begin{defn}\label{tau}
We define a  map $\ta$ from 
the dual Koszul complex to the dual bar complex
with coefficients in $S(V)\# G$:
$$
\Upsilon_p:\ \ 
C^p 
\rightarrow
\Hom_{\CC}(S(V)^{\ot p},S(V)\#G) \ .
$$
For $g$ in $G$ with basis $B_g=\{v_1,\ldots, v_n\}$ of $V$ as above,
and 
$
\alpha = f_g \otimes v_{j_1}^*\wedge\cdots\wedge v_{j_p}^*\ot g$
with 
$f_g\in S(V)$ and $1\leq j_1<\ldots<j_p \leq n$, 
define $\ta(\alpha):S(V)^{\ot p}\rightarrow S(V)\# G$ 
by
$$
\begin{aligned}
\ta(\alpha)(& f_1\otimes \cdots\otimes f_p )
&= \Biggl(\
\prod_{k=1,\ldots, p}\
^{s_{1}s_{2}\cdots s_{j_{k}-1}} 
(\del_{j_{k}}f_k) \Biggr) f_g\ot g
\ .
\end{aligned}
$$
By Theorem \ref{tauproperties} below, $\ta$ is a cochain map. Thus
$\ta$ induces a map on the cohomology
$\HH^{\DOT}(S(V), S(V)\#G)$, which we denote by $\ta$ as well.
\end{defn}

We make the following remark, which will be needed in our analysis
of the cup product in Section~\ref{cup=smash}. 

\vspace{3ex}
\begin{remark}
\label{zeroterm}
For the fixed basis $B_g=\{v_1,\ldots, v_n\}$ 
and $\alpha=f_g \otimes v_{j_1}^*\wedge\cdots\wedge v_{j_p}^*\ot g$ 
in $C^p$ (with $j_1<\ldots< j_p$), note that 
$$
\ta(\alpha)(v_{i_1} \otimes\cdots\ot v_{i_p})=0
\quad\quad\quad
\text{unless } i_1=j_1, \ldots, i_p=j_p\ .
$$
In general, 
$
\ta(\alpha)(f_1\ot\cdots\ot f_p)=0
$
whenever $\displaystyle{\frac{\del}{\del v_{j_k}}(f_k)} =0$ for some $k$. 
\end{remark}
\vspace{3ex}

We shall use the following consequence of the definitions.
  
\vspace{2ex}
\begin{prop}\label{phi*tau=1}
For any choices of bases defining $\ta$,
$$
\Phi^* \ta = 1\ 
$$
as a map on cochains $C^{\DOT}$.
\end{prop}
\vspace{1ex}
\begin{proof}
Consider a nonzero cochain $\alpha$ in $C^p$.
Without loss of generality,
suppose that $\alpha=f_g\otimes v_1^*\wedge\cdots\wedge v_p^*\ot g$
for some $g$ in $G$, where $f_g$ is in $S(V)$ and
$B_g=\{v_1,\ldots, v_n\}$ is the fixed basis of eigenvectors of $g$.
Then for all $1\leq i_1,\ldots, i_p\leq n$,
\begin{flalign*}
&&(\Phi^*\ta)(\alpha)(v_{i_1}\wedge\cdots\wedge v_{i_{p}})
& \quad=\quad
\ta(\alpha)(\Phi(v_{i_1}\wedge\cdots\wedge v_{i_{p}})) &&\\
&&& \quad=\quad
\ta(\alpha)
\ \biggl(\, \sum_{\pi\in \Sym_{p}} \sgn(\pi)\ v_{i_{\pi(1)}}\ot\cdots \ot v_{i_{\pi(p)}}\biggr)\\
&&&\quad=\quad
\sum_{\pi\in \Sym_{p}} \sgn(\pi)\
   \ta(\alpha)(v_{i_{\pi(1)}}\ot\cdots\ot v_{i_{\pi(p)}})\ ,&&\\
\text{while}
&& \alpha(v_{i_1}\wedge\cdots\wedge v_{i_{p}})
&\quad=\quad
f_g \ (v_1^*\wedge\cdots\wedge v_p^*)(v_{i_1}\wedge\cdots\wedge v_{i_{p}})\ot g\ .&&
\end{flalign*}
By Remark~\ref{zeroterm}, both
expressions are zero unless
$\{v_{i_1}, \ldots v_{i_p}\}=\{v_1,\ldots, v_p\}$, in which case
both yield
$\sgn(\pi)\, f_g \ot g$,
where $\pi$ is the permutation sending $(i_1,\ldots, i_p)$ to
$(1,\ldots, p)$.
\end{proof}
\vspace{2ex}

We summarize results needed from~\cite{SWchain}:
\vspace{2ex}
\begin{thm}
\label{tauproperties}
The combinatorial converter map
$$
\begin{aligned}
\ta:\ \ \text{Dual Koszul Complex}
\quad&\rightarrow\quad\text{Dual Bar Complex}\\
C^{p}
\quad&\rightarrow\quad
\Hom_{\CC}(S(V)^{\ot p},S(V)\#G) \ 
\end{aligned}
$$
induces isomorphisms of cohomology independent of choices of bases:
\begin{itemize}
\item
For any  basis $B_g$ of eigenvectors for any $g$ in $G$,
$\ta$ is a cochain map.
\item
Although the cochain map $\ta$ depends on  choices $B_g$
for $g$ in $G$, the induced map $\ta$ on cohomology $\HH^{\DOT}(S(V), S(V)\# G)$ 
is independent of choices.
\item
The map $\ta$ induces an automorphism of $\HH^{\DOT}(S(V),S(V)\#G)$
with inverse automorphism  $\Phi^*$.
Specifically, $\ta$ and $\Phi^*$ convert 
between expressions of cohomology in terms of the Koszul
resolution and the bar resolution.
\item
The map $\ta$ on $\HH^{\DOT}(S(V),S(V)\#G)$ is $G$-invariant
and hence induces an automorphism on
$\HH^{\DOT}(S(V)\#G)\cong \HH^{\DOT}(S(V),S(V)\#G)^G\ .$
\end{itemize}
\end{thm}
\vspace{3ex}

\begin{remark}
We do not symmetrize $\ta$, in comparison with
similar maps in the literature (see Anno~\cite{Anno} and Halbout
and Tang~\cite{HalboutTang}).
Since they are chain maps, 
our maps are the same on
cohomology as their symmetrized versions.
Symmetrization may be more elegant, however unsymmetrized maps can be more convenient
for computation.
\end{remark}
\vspace{3ex}

\section{Smash product structure}\label{smashstructure}

In this section, we transform the
``cohomology of a group action'' into a ``group action on cohomology'' by
viewing both as algebras.  
We relate Hochschild cohomology of the skew group
algebra, $\HH^{\DOT}(S(V)\# G)$, to the skew group algebra
of a Hochschild cohomology algebra, $\HH^{\DOT}(S(V))\# G$:
We manifest the first algebra as a subquotient of the second.
We thus twist the cup product on $\HHD(S(V))$ by the group
action and obtain a natural smash product on $\HHD(S(V)\#G)$.
In the next section, we show that the cup product on
$\HHD(S(V)\#G)$ is precisely this natural smash product.

We first embed $\HH^{\DOT}(S(V), S(V)\#G)$ as a graded vector space
into $$\HH^{\DOT}(S(V))\,\#\, G\ ,$$
the skew group algebra determined by the action of $G$ on 
$\HH^{\DOT}(S(V))$
induced from its action on $V$.  
Note that for a general algebra $S$ with action of $G$ by automorphisms,
$\HH^{\DOT}(S,S\# G)$ is a $G$-{\em graded} algebra;
we show that in the special case $S=S(V)$, it is not only
$G$-graded, but is very close to being a smash product itself.
To see this,  we first identify 
the Hochschild cohomology $\HH^{\DOT}(S(V))$ with the set of vector
forms on $V$ (cf.~\cite{Weibel}):
$$\HH^{\DOT}(S(V))=S(V)\ot \Wedge^{\DOT} V^*\ .$$  
The group $G$ acts on this tensor product diagonally,
and the skew group algebra
$\HH^{\DOT}(S(V))\# G$ is the $\CC$-vector space of cochains,
$$
\HH^{\DOT}(S(V))\ot \CC[G] = S(V)\ot \Wedge^{\DOT} V^* \ot \CC[G] = C^{\DOT}\ ,
$$
together with smash product
\begin{equation}
\label{starproduct}
(f_g \otimes dv_g\ot{g})
\sta (f_h \otimes dv_h\ot{h})
=f_g\, ^g\scooch f_h \otimes (dv_g\wedge\, ^g dv_h) \ot{gh}\ ,\\
\end{equation}
where $f_g, f_h \in S(V)$, $g,h \in G$, 
and $dv_g, dv_h \in \Wedge^{\DOT}V^*$.

We  regard $\HH^{\DOT}(S(V),S(V)\# G)$ as a {\em vector space subquotient}
of $\HH^{\DOT}(S(V))\# G$ by identifying with $Z^{\DOT}/B^{\DOT}$
(see Remark~\ref{H_g} and the comments before it):
\begin{equation}\label{vectorspacesubquotient}
\HH^{\DOT}(S(V),S(V)\# G)= Z^{\DOT}/ B^{\DOT}
\ \ \subset\ \  
C^{\DOT}/B^{\DOT}
=\bigl(\HH^{\DOT}(S(V))\ot \CC[G]\bigr)/B^{\DOT}\ .
\end{equation}
We next recognize this subquotient of vector spaces 
as a {\em subquotient of algebras} under the smash product.

\vspace{2ex}
\begin{prop}\label{closedundersmash}
Under the smash product of $\HH^{\DOT}(S(V))\# G$,
the space $Z^{\DOT}$ forms a subalgebra of $C^{\DOT}$
and the space $B^{\DOT}$ forms an ideal of $Z^{\DOT}$:
For all cochains $\alpha$ and $\beta$ in $\HH^{\DOT}(S(V))\#G$,
\begin{itemize}
\item
If $\alpha$ and $\beta$ lie in $Z^{\DOT}$, then
$\alpha\sta\beta$ also lies in $Z^{\DOT}$,
\item
If $\alpha$ lies in $Z^{\DOT}$ and $\beta$ lies in $B^{\DOT}$, 
then $\alpha\sta\beta$  and $\beta\sta \alpha$ also lie in $B^{\DOT}$.
\end{itemize}
\end{prop}
\vspace{1ex}
\begin{proof}
Let $\alpha$ and $\beta$ be cocycles in $Z^{\DOT}$.
Without loss of generality, suppose that 
$$
\alpha=f_g\ot dv_g \ot \vol_g^\perp\ot g
\quad\text{and}\quad
\beta=f_h\ot dv_h \ot \vol_h^\perp\ot h
$$
for $g,h$ in $G$, 
$f_g, f_h$ in $S(V)$, and $dv_g, dv_h$ in $\Wedge V^*$.
By Lemma~\ref{wedgewithperpspace} (see~(\ref{starproduct})),
$$\begin{aligned}
\alpha\sta\beta
&=f_g\, ^g\scooch f_h \otimes (dv_g\wedge\vol_g^\perp \wedge\, ^g(dv_h\wedge
\vol_h^\perp)) 
\ot {gh}\\
&=f_g\, ^g\scooch f_h \otimes (dv_g\wedge\vol_g^\perp \wedge\, dv_h\wedge
\vol_h^\perp) \ot {gh}\\
&=\pm\ f_g\, ^g\scooch f_h \otimes (dv_g\wedge dv_h \wedge \vol_g^\perp
\wedge\vol_h^\perp) 
\ot {gh}\ .
\end{aligned}
$$
Assume this product is nonzero. 
Then by Proposition~\ref{wedgeofvols} and Lemma~\ref{codims}, 
$$(V^g)^\perp\oplus (V^h)^\perp=(V^{gh})^\perp=(V^{hg})^\perp\ $$
and
$\vol_g^\perp\wedge\vol_h^\perp$ is a scalar multiple of
$\vol_{gh}^\perp$.  Hence
$\alpha\sta\beta$ lies in $Z^{\DOT}$.

Now assume further that $\beta$ is a coboundary in $B^{\DOT}$,
i.e., that $f_h$ lies in the ideal $I((V^h)^\perp)$ of $S(V)$.
Note that
$\, ^g((V^h)^\perp)\subset\, ^g((V^{hg})^\perp)=(V^{ghgg^{-1}})^\perp
=(V^{gh})^\perp$.
Hence $\, ^gf_h$, and thus the product
$f_g \, ^g\scooch f_h$,
 lies in $I((V^{gh})^\perp)$.
Therefore, $\alpha\sta\beta$ is an element of $B^{\DOT}$.
The argument for $\beta\sta\alpha$ is similar (and easier).
\end{proof}
\vspace{2ex}

The proposition above immediately implies that the smash
product on the skew group algebra $\HH^{\DOT}(S(V))\# G$ induces a smash algebra product
on the cohomology $\HH^{\DOT}(S(V), S(V)\#G)$, as we see in the next two results. 

\vspace{2ex}
\begin{cor}\label{Z/Balg}
The vector space subquotient $Z^{\DOT}/B^{\DOT}$ of $\HHD(S(V))\# G$
is an algebra subquotient (subalgebra of a quotient of algebras) under the induced smash product.
\end{cor}
\vspace{2ex}

Note that cohomology classes (and cocycles) in general are represented by 
{\em sums} of elements of the form of 
$\alpha$ and $\beta$ given in the next proposition.

\vspace{2ex}
\begin{prop}\label{inducedsmash}
The cohomology $\HH^{\DOT}(S(V), S(V)\#G)$ identifies naturally as a graded
vector space with an algebra subquotient of the smash
product $\HH^{\DOT}(S(V))\# G$.  Under this identification,
$\HH^{\DOT}(S(V), S(V)\#G)$ inherits the smash product:
For cohomology classes in $\HH^{\DOT}(S(V), S(V)\#G)$ 
represented by cocycles in $Z^{\DOT}$,
$$
\alpha=\sum_{g\in G}f_g \otimes dv_g\ot g\quad\text{and}\quad
\beta=\sum_{h\in G}f_h \otimes dv_h \ot h\ ,
$$ 
where each $f_g, f_h \in S(V)$
and each $dv_g, dv_h \in \Wedge^{\DOT} V^*$,
the smash product
\begin{equation}
\alpha\sta\beta
=\sum_{{g, h \in G}\atop{g\leq gh}}
f_g\, ^g\scooch f_h \otimes (dv_g\wedge\, ^g dv_h)
\ot\ {gh}\ 
\end{equation}
is also a cocycle representing a class of $\HH^{\DOT}(S(V), S(V)\#G)$.
\end{prop}
\vspace{1ex}
\begin{proof}
We saw in~(\ref{vectorspacesubquotient}) that
$\HH^{\DOT}(S(V), S(V)\#G)$ is isomorphic to 
the vector space subquotient 
$$Z^{\DOT}/B^{\DOT}\quad\text{of}\quad \bigl(\HH^{\DOT}(S(V))\# G \bigr)/ B^{\DOT}\ .$$  
(In fact, we may identify
$\HH^{\DOT}(S(V), S(V)\#G)$ with a {\em subset} of 
$\HH^{\DOT}(S(V))\# G$; see
Remarks~\ref{H_g} and~\ref{H_gsmash}.)  By Proposition~\ref{closedundersmash},
the smash product on $\HH^{\DOT}(S(V))\# G$ may be restricted to
the subset $Z^{\DOT}$ and induces a multiplication on
the subquotient $Z^{\DOT}/B^{\DOT}$ (Corollary~\ref{Z/Balg}).  
Hence $\HHD(S(V), S(V)\#G)$ is isomorphic to an algebra carrying
a natural smash product, and thus it inherits (under this isomorphism)
a natural smash product of its own.
Formula (\ref{starproduct}) gives a cohomology class representative
of the smash product of two cocycles.
\end{proof}
\vspace{2ex}

The proposition above explains that the vector space
inclusion of (\ref{vectorspacesubquotient}) yields an injection of algebras:
The algebra $\HH^{\DOT}(S(V),S(V)\# G)$, under smash product~$\sta$,
is isomorphic to an algebra subquotient of $\HH^{\DOT}(S(V))\# G$. 
But for any $\CC$-algebra $A$ carrying an action of $G$, 
the smash product on the skew group algebra $A\#G$ maps $G$-invariants
to $G$-invariants.
Hence, the algebra $\HHD(S(V)\#G)$ also inherits a smash product
and is isomorphic (under smash product) to an algebra subquotient
of $\HHD(S(V))\#G$. We shall see in Corollary~\ref{algebrasubquot} below that the
same is true under cup product.

\vspace{3ex}
\begin{remark}\label{H_gsmash}
One may identify $\HH^{\DOT}(S(V), S(V)\#G)$ with the subset $H^{\DOT}$ of 
the smash product $\HH^{\DOT}(S(V))\# G$ by fixing a set of cohomology class
repesentatives as in Remark~\ref{H_g}.
But note that $H^{\DOT}$ is not closed under the smash product 
as a subset of $\HH^{\DOT}(S(V))\# G$---we 
must take a quotient by coboundaries and again
take chosen representatives for this quotient.
The induced smash product explicitly
becomes
$$
(\alpha,\beta)\mapsto\Proj_H(\alpha\sta\beta)\ .
$$
\end{remark}
\vspace{3ex}


\section{Equivalence of cup and smash products}\label{cup=smash}

In the previous section, we twisted the cup product on $\HHD(S(V))$
by the action of $G$ in a natural way to define a multiplication on 
$\HHD(S(V),S(V)\#G)$ and also on $\HHD(S(V)\#G)$.
In fact, we showed that as a graded vector space,
the Hochschild cohomology $\HH^{\DOT}(S(V)\# G)$ 
maps isomorphically to an algebra subquotient of
$\HH^{\DOT}(S(V))\# G$ 
and thus inherits a natural smash product structure.
We now regard $\HH^{\DOT}(S(V), S(V)\# G)$ as an algebra under three operations:
\begin{itemize}
\item the cup product $\smile$ induced from the bar resolution of $S(V)$, 
\item a ``smash'' product $\sta$ induced from  $\HH^{\DOT}(S(V))\# G$, and
\item the usual multiplication in the algebra tensor product $\HHD(S(V))\ot \CC[G]$.
\end{itemize}
We show in this section that these three basic algebraic operations
coincide.  
This allows us to describe generators for Hochschild cohomology 
in the next section.  The equality of cup and smash products
(given in the theorem below) also explains how the cup product,
defined on the bar resolution, 
may be expressed in terms of the Koszul resolution.

\vspace{2ex}
\begin{thm}\label{thm:cup=smash}
For all cocycles $\alpha$ and $\beta$ in $Z^{\DOT}$,
$$\alpha \smile \beta = \alpha \sta \beta\ .$$
On the Hochschild cohomology $\HH^{\DOT}(S(V), S(V)\# G)$,
this product is induced from the usual
multiplication on the tensor product 
$\HHD(S(V))\ot \CC[G]$ of algebras.
\end{thm}
\vspace{1ex}
\begin{proof}
We need only verify the statement
for cocycles in $Z^{\DOT}$ of the form
$
  \alpha = f_g \ot dv_g\ot{g}$
and
$
  \beta = f_h \ot dv_h\ot{h},
$
for some $g,h$ in $G$, 
where $f_g,f_h \in S(V)$ and $dv_g\in \bigwedge^p V^*$, 
$dv_h\in \bigwedge^q V^*$ (as any cocycle 
in $Z^{\DOT}$ will be the sum of such elements).  

As $\Phi^*$ and $\ta$ are inverse maps on the cohomology
$\HH^{\DOT}(S(V),S(V)\# G)$ converting
between Koszul and bar cochain complexes 
(see Diagram~\ref{big-diagram-tau}),
$$\alpha\smile\beta = \Phi^* (\ta(\alpha)\smile\ta(\beta))$$
where $\smile$ (on the right hand side) denotes cup product on the bar complex.
Suppose
$B_g = \{v_1,\ldots, v_n\}$, a basis of eigenvectors of
$g$ (see Definition~\ref{tau}).
Without loss of generality, assume that 
$dv_g = v_1^*\wedge\cdots\wedge v_p^*$ where the span of
$v_1, \ldots, v_p$ includes $(V^g)^\perp$
(see (\ref{cocycles})).  
(In general,  $dv_g$ will be a sum of such elements
with indices relabeled.)
By the second part of Lemma~\ref{wedgewithperpspace},
$$
\begin{aligned}
\alpha \sta \beta
&=(f_g \otimes dv_g\ot{g})
\sta (f_h\otimes dv_h\ot{h})\\
&=f_g\, ^g\scooch f_h\otimes (dv_g\wedge\, ^g dv_h)\ot{gh}\\
&=f_g\, ^g\scooch f_h\otimes (dv_g\wedge\, dv_h)\ot{gh}\ .
\end{aligned}
$$

We compare the values of $\alpha\smile\beta$ and $\alpha \sta \beta$ 
on any $v_{i_1}\wedge\cdots\wedge v_{i_{p+q}}$. Now
$$
\begin{aligned}
&(\alpha \sta\beta)(v_{i_1}\wedge\cdots\wedge v_{i_{p+q}})
=
f_g\, ^g\scooch f_h\ (dv_g\wedge\, dv_h)
(v_{i_1}\wedge\cdots\wedge v_{i_{p+q}})\ot{gh}\ \\
\text{while}\\
&(\alpha\smile\beta)(v_{i_1}\wedge\cdots\wedge v_{i_{p+q}})\\
& \hspace{5ex}=
\Phi^* (\ta(\alpha)\smile \ta(\beta))
 (v_{i_1}\wedge\cdots\wedge v_{i_{p+q}}) \\
& \hspace{5ex}=(\ta(\alpha)\smile \ta(\beta))
  \bigl(\sum_{\pi\in \Sym_{p+q}} \sgn(\pi) v_{i_{\pi(1)}}\ot\cdots \ot v_{i_{\pi(p+q)}}\bigr)\\
&\hspace{5ex}=\sum_{\pi\in \Sym_{p+q}} \sgn(\pi)\
   \ta(\alpha)(v_{i_{\pi(1)}}\ot\cdots\ot v_{i_{\pi(p)}})\ 
   \ta(\beta)(v_{i_{\pi(p+1)}}\ot\cdots\ot v_{i_{\pi(p+q)}})\ .
\end{aligned}
$$
By Remark~\ref{zeroterm},
both $\alpha\smile\beta$ and $\alpha \sta\beta$ are
readily seen to be zero on $v_{i_1}\wedge\cdots\wedge v_{i_{p+q}}$
unless $\{v_{i_1},\ldots, v_{i_{p+q}}\}$ contains
$\{v_1,\ldots, v_p\}$. 
Thus, we may assume (after relabeling indices and also 
possibly changing signs throughout)
that $v_1=v_{i_1}, \ldots, v_{p+q}=v_{i_{p+q}}$.

Again by Remark~\ref{zeroterm}, $(\alpha\smile\beta)
(v_1\wedge\cdots\wedge v_{p+q})$ is equal to
$$
\sum_{\pi\in \Sym_p}
\sgn(\pi)\
\ta(\alpha)(v_{\pi(1)}\otimes\cdots\otimes v_{\pi(p)})
\sum_{\pi\in \Sym_q}
\sgn(\pi)\
\ta(\beta)(v_{p+\pi(1)}
\ot\cdots\ot v_{p+\pi(q)})\ .
$$
Proposition~\ref{phi*tau=1} then implies that this value of the cup product is just
$$
\begin{aligned}
\Phi^* (\ta(\alpha))
& (v_1\wedge\cdots\wedge v_{p})
\cdot
\Phi^* (\ta(\beta))
 (v_{p+1}\wedge\cdots\wedge v_{p+q})\\
&\hspace{.5cm}=
\alpha
(v_1\wedge\cdots\wedge v_{p})\cdot
\beta
(v_{p+1}\wedge\cdots\wedge v_{p+q})\\
&\hspace{.5cm}=
f_g\, ^g\scooch f_h\ 
dv_g(v_1\wedge\cdots\wedge v_{p})\
dv_h(v_{p+1}\wedge\cdots\wedge v_{p+q})\ot gh\\
&\hspace{.5cm}=
f_g\, ^g\scooch f_h\ (dv_g\wedge\, dv_h)
(v_{1}\wedge\cdots\wedge v_{p+q})\ot gh\\
&\hspace{.5cm}=
(\alpha \sta\beta)(v_1\wedge\cdots\wedge v_{p+q})\ 
\end{aligned}
$$
as elements of $S(V)\#G$.
Thus $\alpha\smile\beta$ and $\alpha\sta\beta$ agree as cochains.

If this product is nonzero, then
$f_h-\, ^g f_h$ lies in the ideal $I((V^g)^{\perp})
\subset I((V^{gh})^{\perp})$
by Lemma~\ref{codims}, Proposition~\ref{wedgeofvols}, and 
the first part of Lemma~\ref{wedgewithperpspace}.
Hence,
$$
f_g\, ^g\scooch f_h\ \ot (dv_g\wedge\, dv_h)
\ot gh \quad\text{and}\quad
f_g\, f_h\ \ot(dv_g\wedge\, dv_h)\ot gh$$
represent the same class in the cohomology
$\HHD(S(V),S(V)\#G)=Z^{\DOT}/B^{\DOT}$.
Thus the usual multiplication 
in the tensor product of algebras $\HHD(S(V))\ot\CC[G]$
gives the cup product on cohomology.
\end{proof}
\vspace{1ex}

We next give an example to show that the cup and smash products do not
agree on arbitrary cochains.

\vspace{3ex}
\begin{example}\label{elem-abel}
Let $G=\ZZ/2\ZZ \times \ZZ/2\ZZ\times \ZZ/2\ZZ$, generated by elements $a_1,a_2,a_3$. Let $V=\CC^3$
with basis $v_1,v_2,v_3$ on which $G$ acts as follows:
$$
  {}^{a_i}v_j = (-1)^{\delta_{ij}}v_j.
$$
Let $\alpha=1\ot v_3^*\ot a_1a_3$ and $\beta = 1\ot v_1^*\wedge v_2^*\ot
a_1a_2$
in $C^{\DOT}$.
Then $$\alpha\sta \beta = -1\ot v_1^*\wedge v_2^*\wedge v_3^*\ot a_2a_3=
-\alpha\smile\beta.$$
\end{example}
\vspace{3ex}

Our results imply the following explicit formula for the cup product, first given by 
Anno~\cite{Anno}, expressed here in terms of the poset in Definition~\ref{binaryrelation}
and the multiplicative cocycle
$\vartheta:G\times G\rightarrow \CC$
of Proposition~\ref{wedgeofvols}.
Note that the condition $g<gh$ in the sum below is included merely
for computational convenience, as 
$\vartheta(g,h)$ is nonzero if and only if $g<gh$. 
Also note that cohomology classes (and cocycles) in general are represented by 
{\em sums} of such $\alpha$ and $\beta$ given in the theorem
below.  

\vspace{2ex}
\begin{thm}\label{cupformula}
Consider cohomology classes in $\HH^{\DOT}(S(V), S(V)\#G)$
represented by
$$
\alpha=\sum_{g\in G}f_g \otimes dv_g\ot \vol_g^\perp\ot g\quad\text{and}\quad
\beta=\sum_{h\in G}f_h \otimes dv_h\ot \vol_h^\perp\ot h\ 
$$ 
in $Z^{\DOT}$ 
of degrees $p$ and $q$, resp.\ (see~(\ref{cocycles})).
Set $m=(\codim V^g)(q-\codim V^h)$.
The cup product $\alpha\smile\beta$
is represented by the cocycle
\begin{equation}
\label{eqn:starproduct}
\sum_{{g, h \in G}\atop{\rule[0ex]{0ex}{1ex}
\text{with} \ g\leq gh}}(-1)^m\, \vartheta(g,h)\
f_g\, f_h \otimes (dv_g\wedge dv_h)\ot \vol^{\perp}_{gh} 
\ot\ {gh}\ .
\end{equation}
\end{thm}
\vspace{1ex}
\begin{proof}
The formula follows directly from 
Proposition~\ref{wedgeofvols}, Proposition~\ref{inducedsmash},
and Theorem~\ref{thm:cup=smash} (see its proof). 
Note that the factor $(-1)^m$ arises when we replace $\vol^{\perp}_g\wedge
\,dv_h$ with $dv_h \wedge\vol^{\perp}_g$.
\end{proof}

\vspace{3ex}
\begin{example}
Let $G,V$ be as in Example~\ref{elem-abel}.
Let $\alpha = 1\ot v_1^*\ot a_1$, $\beta=1\ot v_3^*\wedge v_2^*\ot a_2$
in $H^{\DOT}$.
Let $\vol^{\perp}_{a_1} = v_1^*$, $\vol^{\perp}_{a_2}=v_2^*$, and
$\vol^{\perp}_{a_1a_2} = v_1^*\wedge v_2^*$, so that $\vartheta(a_1,a_2)=1$.
Then 
\begin{eqnarray*}
  \alpha\smile\beta \ \ = \ \  \alpha\sta\beta & = & 
   (-1)^{1(2-1)} \ot  {}^{a_1}(v_3^*)\wedge v_1^*\wedge v_2^*\ot a_1a_2 \\
    & = & -1\ot v_3^*\wedge \vol^{\perp}_{a_1a_2}\ot\ a_1a_2.
\end{eqnarray*}
\end{example}
\vspace{3ex}

\vspace{2ex}
\begin{cor}\label{algebrasubquot}
The algebra $\HH^{\DOT}(S(V),S(V)\# G)$ is isomorphic to an 
algebra subquotient of $\HH^{\DOT}(S(V))\# G$. 
The algebra $\HH^{\DOT}(S(V)\# G)$ is isomorphic to its $G$-invariant
subalgebra.  Hence, it is 
also an algebra subquotient of $\HHD(S(V))\#G $.
\end{cor}
\vspace{1ex}
\begin{proof}
The space $\HHD(S(V),S(V)\#G)$ may be written as a (graded) subspace
of a vector space quotient of $\HHD(S(V))\#G$ (see~(\ref{vectorspacesubquotient})). 
Proposition~\ref{inducedsmash} states that this subquotient
is actually a subquotient of algebras under the smash product
(see Proposition~\ref{closedundersmash} and Corollary~\ref{Z/Balg}).
But the cup product on $\HHD(S(V),S(V)\# G)$ is the same as this
induced smash product by Theorem~\ref{thm:cup=smash}.
By Theorem~\ref{isos}, 
the cohomology algebra $\HH^{\DOT}(S(V)\# G)$ is isomorphic to
the $G$-invariant subalgebra of $\HH^{\DOT}(S(V),S(V)\# G)$ (both
spaces regarded as algebras under their respective cup products).
Hence the algebra $\HHD(S(V)\#G)$ is also a subquotient of $\HHD(S(V))\#G$.
\end{proof}
\vspace{3ex}

\section{Volume algebra}\label{volumealgebra}
Results in the last section reveal an interesting subalgebra
of the Hochschild cohomology $\HH^{\DOT}(S(V),S(V)\#G)$ isomorphic
to algebras 
that have appeared in the literature before (see Remark~\ref{describesubalg} below);
 we call it the {\bf volume algebra}.  
We use this subalgebra to give algebra generators of
cohomology in the next section.

We define the volume algebra with the next proposition.
Recall that for each $g$ in $G$, the form $\vol^{\perp}_g$ is a choice of nonzero
element in the one-dimensional space $\Wedge^{\codim V^g}((V^g)^{\perp})^*$.
We show that the elements
$\vol_g^\perp\ot\, g:= 1 \ot \vol_g^\perp\ot\, g$ 
in $\HHD(S(V),S(V)\#G)$ generate a subalgebra
which captures the binary relation $\leq$ on $G$ and 
reflects the poset structure of $G/K$
(see Definition~\ref{binaryrelation}).
(Note that if $k$ acts trivially on $V$,
then $\vol_k^{\perp}=1$ up to a nonzero constant.)

\vspace{2ex}
\begin{prop}\label{subalg}
The $\CC$-vector space
$$
A_{\vol} := \Span_{\CC}
\{\vol_g^\perp\ot\, {g} \mid g\in G\}\ 
$$
is a subalgebra of $\HH^{\DOT}(S(V),S(V)\# G)$.
The induced multiplication on this subalgebra is given by
$$
(\vol_g^\perp\ot\,{g})(\vol_h^\perp\ot\, h)
=
\vartheta(g,h)
\vol_{gh}^\perp\ot\,{gh}
$$ 
where 
$\vartheta: G\times G\rightarrow \CC$ is a (multiplicative) cocycle on $G$.
For any $g,h$ in $G$, the above product 
is nonzero ($\vartheta(g,h)\neq 0$)
if and only if
$g\leq gh$.
The algebra $A_{\vol}$
is isomorphic to the (generalized) twisted group 
algebra $\CC^{\vartheta}[G]$.
\end{prop}
\vspace{1ex}
\begin{proof}
Lemma~\ref{wedgewithperpspace} and  
Proposition \ref{wedgeofvols} 
imply that for any pair $g,h$ in $G$,
$$
\begin{aligned}
(\vol_g^\perp\ot\, g)
\ \sta \
(\vol_h^\perp\ot\, h)
& =
(\vol_g^\perp\wedge\ ^g(\vol_h^\perp))\ot\, gh \\
& =
(\vol_g^\perp\wedge\, \vol_h^\perp)\ot\, gh \\
&=  \vartheta(g,h) \vol_{gh}^\perp\ot\, gh\ ,
\end{aligned}
$$
for a cocycle $\vartheta$ which is nonzero on the pair $(g,h)$
exactly when $g\leq gh$. 
Hence by Theorem~\ref{thm:cup=smash},
$A_{\vol}$ is closed under cup product
and isomorphic to $\CC^{\vartheta}[G]$.
\end{proof}
\vspace{2ex}

The above proposition explains how 
$A_{\vol}$  algebraically captures  
the geometric relation on $G$ given by $\leq$.
Indeed, the space $A_{\vol}$ forms a {\em graded algebra} where
$$\deg(\vol_g^\perp\ot\, g)=\codim V^g\quad\text{for each}\ g \in G\ .$$
Moreover, for all $g,h$ in $G$,
$$
g\leq h \quad\ \ \text{if and only if}\quad \ \  
\vol_g^\perp \ot\, g \ \text{ divides}\ 
\vol_h^\perp\ot\, h \quad\text{in}\  A_{\vol}\ .
$$
Let $S$ be the subset 
$\{\vol_g^\perp\ot\, g:g \in G\}$ of $A_{\vol}$.
Then $S$ is a poset under division:
$$
\vol_g^\perp \ot\, g \ \leq_{S}\ \vol_h^\perp\ot\, h 
\ \ \quad\text{if and only if}\quad \ \  
\vol_g^\perp \ot\, g \ \ \text{divides}\ \ 
\vol_h^\perp\ot\, h \quad\text{in}\ \ A_{\vol}\ .
$$
If $G$ acts faithfully, we have isomorphic posets:
$$
(G,\leq) \cong (S, \leq_S)\ .
$$

\vspace{3ex}
\begin{remark}\label{describesubalg}
The subalgebra $A_{\vol}$ of $\HH^{\DOT}(S(V), S(V)\# G)$
appears in other arenas.
If $V$ is a symplectic vector space 
and $G\subset {\rm Sp}(V)$, 
the algebra $A_{\vol}$ is isomorphic to the graded algebra
${\rm gr}_F \CC[G]$ associated to the
filtration $F$ on $G$ assigning to each group
element $g$ the codimension of $V^g$.
Its $G$-invariant subalgebra $A_{\vol}^G$ 
is isomorphic to 
the orbifold cohomology of $V/G$ in this case
(see Fantechi and G\"ottsche~\cite[\S2]{Fantechi-Gottsche}
or Ginzburg-Kaledin~\cite{GinzburgKaledin})
and also to
the Hochschild cohomology of the $G$-invariant subalgebra of
the Weyl algebra $A(V)$
(see Suarez-Alvarez~\cite{Suarez-Alvarez}):
$$
\HH^{\DOT}( A(V)^G)\cong \HH^{\DOT}\bigl(A(V)\# G\bigr)
\ \cong \
A_{\vol}^G
\ \cong\
\coh^{\DOT}_{\text{orb}}\bigl(V/G\bigr)\ .
$$
Work of several authors has shown that for $G=\Sym_n$
acting symplectically on $\CC^{2n}$,
this orbifold cohomology is isomorphic
to the cohomology of a Hilbert scheme which is a crepant
resolution of the orbifold
\cite{Fantechi-Gottsche,lehn-sorger01,uribe}.
Lehn and Sorger gave a description of $A^{\Sym_n}_{\vol}$
in terms of generators and relations 
\cite[Remark 6.3]{lehn-sorger01}.
\end{remark}
\vspace{3ex}

\section{Cohomology generators given by the poset}\label{poset}

In this section, we give generators for the algebra $\HH^{\DOT}(S(V), S(V)\# G)$
in terms of the partial order $\leq$ on the quotient group $G/K$, where $K$ is the kernel
of the representation of the group $G$ acting on $V$.
In Section~\ref{cup=smash}, we showed that the cup product and an induced
smash product on the algebra $\HH^{\DOT}(S(V), S(V)\# G)$ agree.
Hence, we simply discuss generation of cohomology as an algebra,
without explicit reference to the product.  We explain in the next three
results how generators for $\HH^{\DOT}(S(V),S(V)\#G)$
are tagged by $K$ together with minimal elements in the poset $G/K$.  
Actually, $1_{G/K}:=K$ is the unique minimal element in the poset $G/K$, 
and we remove it before seeking minimal elements.

\vspace{2ex}
\begin{thm}\label{generatersofA}
The subalgebra $A_{\vol}$ of $\HH^{\DOT}(S(V),S(V)\# G)$ is generated 
as an algebra over $\CC$ by
$$\{\vol_g^\perp\ot\,{g}\mid g\in K \mbox{ or } 
g \in [G/ K]  \mbox{ with $gK$ minimal in the poset } G/ K-\{1_{G/K} \}\} \ ,$$
where $[G/ K]$ is any set of coset representatives for $G/ K$.
\end {thm}
\vspace{1ex}
\begin{proof}
Assume $h\notin K$.  We first observe that we may write $\vol^{\perp}_h\ot\, h$
as the product of two volume forms, one tagged by any other coset representive
of $hK$ and the other tagged by an element of $K$:
Suppose $hK=h'K$ for some 
$h'$ in $G$. Then $V^h=V^{h'}$ and $\vol_{h}^\perp=\vol_{h'}^{\perp}$
up to a nonzero scalar.  Let $k=h^{-1}h'$, so $\vol_k^\perp$ is a nonzero
scalar itself.
Then (see Equation~\ref{starproduct}),
$$
(\vol^{\perp}_h\ot\, h)\sta(\vol^{\perp}_k\ot\, k)
= (\vol_h^{\perp}\wedge\ ^h(\vol_k^\perp))\ot\, hk,
$$
which is a nonzero scalar multiple of 
$\vol_h^{\perp}\ot\, h'$, and thus of $\vol_{h'}^{\perp}\ot\, h'$.
We may therefore assume that $h$ itself lies in $[G/ K]$.

Suppose that $hK$ is not minimal in the partial order on 
$G/ K-\{1_{G/K}\}$.  
Then there exists $g\in [G/K]$ with $K\neq gK\neq hK$
and $gK<hK$, i.e., $(V^g)^\perp\oplus (V^{g^{-1}h})^\perp = (V^h)^\perp$.
By Theorem~\ref{cupformula},
$$ 
(\vol_g^\perp\ot\, g)
\ \sta \
(\vol_{g^{-1}h}^\perp\ot\, g^{-1}h)
= \vartheta(g,g^{-1}h)
\vol_{h}^\perp\ot\, h\ ,
$$
which is nonzero since $gK < hK$.
By symmetry in the definition of the partial order,  
$g^{-1}hK<hK$ in the poset $G/ K-\{1_{G/K}\}$ as well.
Hence, we have written $\vol^{\perp}_h\ot h$
as a product of volume forms each tagged by a group element 
less than $h$ in the partial order.
As the set ${G /K}-\{1_{G/K}\}$ is finite, the partial order is well-founded,
i.e., every descending chain contains a minimal element.
Hence, by induction, $\vol_h^\perp\ot h$ may be written (up to a scalar) 
as the product of elements in the set given in the statement of the theorem.
\end{proof}
\vspace{2ex}

We now turn to the task of describing generators for the full Hochschild
cohomology 
$\HHD(S(V),S(V)\# G)$ as an algebra.
We regard $\HH^{\DOT}(S(V))$ as a subalgebra of
$\HHD(S(V),S(V)\# G)$ by identifying 
$$\HH^{\DOT}(S(V))=S(V)\ot\Wedge^{\DOT}V^*
\quad\text{with}\quad
S(V)\ot\Wedge^{\DOT}V^*\ot 1_G\ .$$ 
Then $\HH^{\DOT}(S(V),S(V)\# G)$ becomes a module over $\HH^{\DOT}(S(V))$ under
cup product. 

\vspace{2ex}
\begin{thm}\label{generatedbysubalg}
The Hochschild cohomology algebra
$\HH^{\DOT}(S(V),S(V)\# G)$ is generated 
by its subalgebras $\HH^{\DOT}(S(V))$ and $A_{\vol}$.
\end{thm}
\vspace{1ex}
\begin{proof}
Let $\alpha$ be an arbitrary
element of $\HH^{\DOT}(S(V),S(V)\# G)$.  
Without loss
of generality, assume
$\alpha = f_h\otimes dy\wedge \vol_h^\perp\ot\, h$ for some
$h$ in $G$ and $dy\in \Wedge ^{\DOT}(V^h)^*$.
By Theorem~\ref{thm:cup=smash} (see~(\ref{eqn:starproduct})), 
$$\alpha = (f_h\ot dy\ot 1_G)\sta (\vol_h^\perp\ot\, h)$$ 
where $f_h\ot dy\ot 1_G$ identifies with $f_h\ot dy$ in $\HH^{\DOT}(S(V))$.
Hence, $\alpha$ lies in the product
$\HH^{\DOT}(S(V))\cdot A_{\vol}$.
\end{proof}
\vspace{1ex}

Recall that for $k$ in $K$, 
$\vol_k^{\perp}=1$ up to a nonzero constant in $\CC$.
The last two theorems then imply:

\vspace{2ex}
\begin{cor}\label{generatednonfaithful}
The Hochschild cohomology algebra
$\HH^{\DOT}(S(V),S(V)\# G)$ is generated 
by its subalgebra $\HH^{\DOT}(S(V))$ and
$$\{\vol_g^\perp\ot\, g\mid g\in K \mbox{ or } 
g \in [G/ K]  \mbox{ with $gK$ minimal in the poset } G/ K-\{1_{G/K}\} \}\ ,$$
where $[G/ K]$ is any set of coset representatives of $G/ K$.
\end{cor}
\vspace{2ex}

If $G$ acts faithfully on $V$, then we may simply take minimal elements
in the poset $G-\{1_G\}$ (see Lemma~\ref{posetlemma}) 
to obtain a generating set under cup product: 
\vspace{2ex}
\begin{cor}\label{generatedfaithful}
Assume  $G$ acts faithfully on $V$.
The Hochschild cohomology algebra
$\HH^{\DOT}(S(V),S(V)\# G)$ is generated 
by $\HH^{\DOT}(S(V))$ and
$$\{\vol_g^\perp\ot\, g\mid  
g \mbox{ is minimal in the poset } G-\{1_G\} \}.$$
\end{cor}
\vspace{3ex}

\section{Reflection groups}\label{reflectiongroups}

In the previous sections, we described the
algebraic structure of the Hochschild cohomology
$\HHD(S(V),S(V)\#G)$.  These results have a special interpretation for reflection groups and
Coxeter groups in particular.
We are interested in comparing the codimension of the fixed point space $V^g$
of a group element $g$ with its ``reflection length'' in the group.

Recall that a nonidentity element of $\GL(V)$ is a reflection
if it fixes a hyperplane in $V$ pointwise.  
A {\bf reflection group} is a finite group generated by reflections.
A reflection group is called a 
{\bf Coxeter group} when it is 
generated by reflections acting on a {\em real} vector space.  
In this section, we restrict ourselves to the case
when $G$ is a reflection group.
We define a length function with respect
to the set of {\em all} reflections inside $G$.
 (Note that this definition may differ from the length function defined
in terms of a fixed choice of generators for the group $G$, for example, a choice
of simple reflections for a Weyl group.)

\vspace{3ex}
\begin{defn}
For each $g$ in $G$, let $l(g)$ be the minimal number $k$ 
such that $g=s_1\cdots s_k$ for
some reflections $s_i$ in $G$.  
We set $l(1_G)=0$. We call $l:G\rightarrow\NN$ the {\bf reflection
  length function} (or ``absolute length function'') 
of $G$.
\end{defn}
\vspace{3ex}

The reflection length function induces a partial order $\leq_l$ on $G$:
$$g\leq_l h \ \ \text{ when }\ \  l(g)+l(g^{-1}h)=l(h)\ .$$  
The poset formed by the reflection length function
plays an important role in the emerging
theory of Artin groups of finite type.  This theory relies on a key result
for Coxeter groups asserting that the closed interval from the identity group element to a
Coxeter element forms a lattice.  
Brady and Watt~\cite{BradyWatt2} gave a case-free proof of
this fact by relating the two partial orders $\leq$ and $\leq_l$.
The poset $\leq_l$ defined from the reflection length function has received attention not
only for Coxeter groups, but for other complex reflection groups $G$ as well.
One may define a Coxeter element and again consider the interval (often a lattice)
from the identity to the Coxeter element, the so-called 
{\em poset of noncrossing partitions for $G$}.
See, for example, Bessis and Reiner~\cite{BessisReiner}.

Note that the length of a linear transformation
with respect to the ambient group $\GL(V)$
coincides with the codimension of the fixed point space:
Any element $g$ in the unitary group $\U(V)$ can be written as the product of 
$m$ reflections in $\GL(V)$ and no fewer if and only if $\codim V^g=m$ 
(see Brady and Watt~\cite{BradyWatt}).

 Reflection length in the group $G$ is bounded below by the codimension of the fixed point
space:
\begin{equation}
l(g) \geq \codim V^g\quad \text{for all}\ \  g \text{ in } G\ .
\end{equation}
Indeed, if $l(g)=m$ and $g=s_1\cdots s_m$ is a product of reflections in $G$,
then
$$V^g\supset V^{s_1}\cap\cdots \cap V^{s_m}\ ;$$
but each $V^{s_i}$ is a hyperplane, so
$$\codim V^g \leq \codim (V^{s_1}\cap\cdots\cap V^{s_m})\leq m\ .$$

For certain reflection groups, reflection length coincides with
codimension of fixed point space. 
The arguments of Carter~\cite{Carter} for Weyl groups hold for
Coxeter groups as well:
\vspace{2ex}
\begin{lemma}\label{carterlemma}
Let $G$ be a (finite) Coxeter group.  Then reflection length coincides with
codimension of fixed point space: For all $g$ in $G$,
$$
\codim V^g = l(g)\ .
$$
\end{lemma}
\vspace{2ex} 

The lemma above implies that 
for Coxeter groups (which act faithfully), 
the partial order $\leq$ above (see Definition~\ref{binaryrelation}) describing 
the ring structure of Hochschild cohomology
coincides with the partial order $\leq_l$ induced from the reflection length
function.
When these two posets agree (i.e., when $\leq\ =\ \leq_l$), 
we may express the ring structure of Hochschild
cohomology in an elegant way.
In particular, the Hochschild
cohomology algebra $\HHD(S(V),S(V)\# G)$ is generated 
in degrees~$0$ and~$1$, as we see in the next corollary.
Note the analogy with the Hochschild-Kostant-Rosenberg Theorem,
which implies that the Hochschild cohomology of a smooth commutative
algebra is also generated in degrees~$0$ and~$1$.

\vspace{2ex}
\begin{thm}\label{degrees0and1}
Suppose $G$ is a (finite) reflection group for which
the reflection length function gives codimension of fixed point spaces:
$$
l(g)=\codim V^g\quad\text{for all }\ g \text{ in } G\ .
$$
Then
$\HH^{\DOT}(S(V),S(V)\# G)$ is generated (as an algebra)
in degrees~$0$ and~$1$.
\end{thm}
\vspace{1ex}
\begin{proof}
As $G$ is a reflection group, it acts faithfully on $V$ by definition.
The elements of $G-\{1_G\}$ minimal in the partial order $\leq_l$
induced by the length function are the reflections.
Indeed, suppose the length of $h$ in $G$ is $m>0$ and write
$h=s_1\cdots s_m$ for some reflections $s_i$ in $G$.
As $h$ can not be expressed as the product of fewer than $m$ reflections,
$l(s_2\cdots s_m)=m-1$.
Hence, 
$$l(s_1)+l(s_1^{-1}h)=l(s_1) + l(s_2\cdots s_m)=1+(m-1)=l(h)$$
and $s_1\leq_l h$.
Note that for any reflection $s$, the relation $g\leq_l s$ implies
that either $g=s$ or $g=1_G$.

By hypothesis, 
the length function and codimension function induce the same partial order.
Hence, the reflections are precisely the minimal elements of 
$G-\{1_G\}$  in the partial order 
$\leq$. 
By Corollary~\ref{generatedfaithful}, 
$\HHD(S(V),S(V)\# G)$ is generated by $\HHD(S(V))$ and
by all $\vol_s^\perp\ot\, s$ where $s$ is a reflection in $G$.
The elements  $\vol_s^\perp\ot\, s$  each have cohomological degree 1, and
$\HHD(S(V))=S(V)\ot \Wedge^{\DOT} V^*$ is generated 
as an algebra (under cup product) by
$\HH^0(S(V))\cong S(V)$ and $\HH^1(S(V))\cong S(V)\ot V^*$.
The statement follows.
\end{proof}

The above corollary applies not only to Coxeter groups, but to other complex
reflection groups as well.  Let $G(r,1,n)$ be the infinite family of complex reflection
groups, each of which is the symmetry group of 
a regular (``Platonic'') polytope in complex space $V=\CC^n$.
The group $G(r,1,n)$ consists of all those $n$ by $n$ complex matrices 
which have in each row and column a single nonzero entry, 
necessarily a primitive $r$-th root of unity.
This group is a natural wreath product
of the symmetric group and the cyclic group of order $r$: $G\cong \ZZ/r\ZZ \wr
\Sym_n$.  In fact, $G(1,1,n)$ is the symmetric group $\Sym_n$
and $G(2,1,n)$ is the Weyl group of type $B_n$.
\vspace{2ex}
\begin{lemma}\label{G(r,1,n)}
For the infinite family $G(r,1,n)$,
the reflection length function coincides with codimension of fixed point spaces:
$$
l(g)=\codim V^g\quad\text{for all}\quad g \text{ in } G\ .
$$
\end{lemma}
\vspace{1ex}
\begin{proof} 
Let $\xi$ be a primitive $r$-th
root of unity in $\CC$.   
 Every element in $G(r,1,n)$ is conjugate to a product $h=c_1\cdots c_k$ 
of disjoint cycles $c_k$ of the form
$$c_k=\xi_j^a\, (i, i+1, \ldots, j)$$
(i.e., $h$ is block diagonal, with $k$-th block $c_k$)
where $\xi_j:=\diag(1,\ldots, 1, \xi, 1,\ldots, 1)$ is the diagonal reflection
with $\xi$ in the $j$-th entry, $a\geq 0$, and where $(i, i+1, \ldots, j)$ is the matrix 
(in the natural reflection representation) of
the corresponding cycle in $\Sym_n$.  (See, for example, Section 2B of
Ram and the first author~\cite{RamShepler}.)
Consider a fixed cycle $c=c_k$ as above.
We may write the cycle $(i, i+1, \ldots, j)$ as a product of $j-i$
transpositions (reflections) in $\Sym_n$, a subgroup of $G(r,1,n)$.
Hence if $\xi^a=1$, then 
$c$ may be expressed as the product of $j-i$ 
reflections in $G(r,1,n)$ while $\codim V^c=j-i$.
If $\xi^a \neq 1 $, then $c$ may be expressed
as the product of $j-i+1$ reflections in $G(r,1,n)$
while $\codim V^c=j-i+1$.  In either case, $l(c)\leq \codim V^c$.
But reflection length is bounded below by codimension,
$\codim V^c \leq l(c)$, and hence $l(c)=\codim V^c$. Then
$$
l(h)\leq l(c_1)+\ldots+ l(c_k)
=\codim V^{c_1}+\ldots + \codim V^{c_k}
=\codim V^h
\leq l(h)
$$
and thus $l(h)=\codim V^h$.
\end{proof}
\vspace{1ex}

As Theorem~\ref{degrees0and1} applies to Coxeter groups by
Lemma~\ref{carterlemma} and to the infinite family $G(r,1,n)$ 
by Lemma~\ref{G(r,1,n)}, we have the following analog of
the Hochschild-Kostant-Rosenberg Theorem:

\vspace{2ex}
\begin{cor}\label{generatorsforCoxeter}
Let $G$ be a Coxeter group or the infinite family $G(r,1,n)$.
Then $\HH^{\DOT}(S(V),S(V)\# G)$ is generated as an algebra in degrees~$0$ and~$1$.
\end{cor}
\vspace{2ex}

Note that the two partial orders $\leq$ and $\leq_l$
do not always agree, i.e., that for a complex reflection group $G$ in general,
the reflection length function may not
give codimension of fixed point spaces:

\vspace{3ex}
\begin{example}\label{Victor}
Let $G$ be the complex reflection group $G(4,2,2)$, the subgroup of $G(4,1,2)$
consisting of those matrices with determinant $\pm 1$.   Let $g$ be the
diagonal matrix $\diag(i,i)$ where $i=\sqrt{-1}$ with determinant $-1$.  
Every reflection in $G$
has determinant $-1$, and hence $g$ can only be written as the product of an
odd number of reflections.  Then $\codim V^g=2$ and yet 
$g$ can not be written as the product of two reflections. 
\end{example}
\vspace{3ex}

\section{Cup product on the invariant subalgebra}\label{last-section}

In the above sections, we investigated the cup product on
the Hochschild cohomology $\HH^{\DOT}(S(V),S(V)\#G)$.  In this section, we describe
the cup product on $\HH^{\DOT}(S(V)\# G)$, its $G$-invariant
subalgebra, using standard techniques from group theory.  
Note that generators for the algebra $\HHD(S(V),S(V)\#G)$ may not be invariant
under $G$ and hence do not generally yield generators
for $\HHD(S(V)\#G)$.

The cup product on $\HH^{\DOT}(S(V)\# G)$ in cohomological degrees $0$ and $1$
is easy to describe, and follows from an observation of
Farinati~\cite{Farinati} (see also~\cite{SheplerWitherspoon}):
\vspace{2ex}
\begin{lemma}\label{lem:determinant}
The only group elements in $G$ that contribute to the Hochschild
cohomology $\HHD(S(V)\#G)$ are those which act on $V$ with determinant $1$.
\end{lemma}
\vspace{2ex}
\noindent
Note that if $G$ embeds in ${\rm SL}(V)$, then
every $g$-component is  nonzero (as any element in 
the one-dimensional subspace $1\ot\Wedge^{\dim V} V^*\ot g$ of $C^{\DOT}$
is automatically invariant under the centralizer of $g$; see~(\ref{eqn:alternate}) below).

Lemma~\ref{lem:determinant} together with Remark~\ref{H_g}
immediately implies that the cup product 
on $\HH^{\DOT}(S(V)\#G)$ in cohomological degrees $0$ and $1$ is
simply the exterior (wedge) product of forms when $G$ acts faithfully
(since reflections do not have determinant~1):

\vspace{2ex}
\begin{prop}
Assume $G$ acts faithfully on $V$, that is, $G$ embeds in
${\rm GL}(V)$. Then the cup product of elements in $\HH^0(S(V)\#G)$
and $\HH^1(S(V)\#G)$ 
is given by the cup product on the $G$-invariant subspaces
of $\HH^0(S(V))$ and $\HH^1(S(V))$:
$$\HH^0(S(V)\# G)=S(V)^G\quad\text{and}\quad
\HH^1(S(V)\#G)=(S(V)\otimes V^*)^G\ .$$
\end{prop}
\vspace{2ex}

When $G$ acts nonfaithfully, we similarly find that only the kernel
$K$ of $G$ acting on $V$ contributes
to the cohomology $\HHD(S(V)\#G)$ in degrees~$0$ and~$1$:
$$\HH^0(S(V)\#G)=\biggl(\,\bigoplus_{k\in K}S(V)\ot k\biggr)^G
\ \ \text{and}\ \
\HH^1(S(V)\# G)=\biggl(\,\bigoplus_{k\in K} S(V)\ot V^*\ot k\biggr)^G\ .$$

The cup product in higher degrees is not as transparent.  We give a formula
in terms of a fixed set 
$\mathcal C$ of representatives of the
conjugacy classes of $G$.
We extend the isomorphism of 
Theorem~\ref{isos}, $\HHD(S(V)\# G)\cong \HHD(S(V),S(V)\# G)^G$:
\begin{equation}\label{eqn:alternate}
  \HHD(S(V)\# G)\cong\bigoplus_{g\in {\mathcal C}} \HHD(S(V),S(V)
  \otimes g)^{Z(g)},
\end{equation}
where $Z(g)$ denotes the centralizer in $G$ of $g$.
The term
$\HHD(S(V),S(V)\otimes g)$ is isomorphic to the $g$-component
of $H^{\DOT}$ of Remark~\ref{H_g}.
The isomorphism 
identifies an element $\alpha$ of $\HHD(S(V),S(V)\otimes g)^{Z(g)}$ 
with the sum $$\sum_{h\in [G/Z(g)]}{}^h \alpha,$$ 
where $[G/Z(g)]$ denotes a set of representatives of left cosets of
$Z(g)$ in $G$.
In the next proposition, we give a formula 
for the cup product of $\HHD(S(V)\#G)$ expressed in terms of the additive decomposition~(\ref{eqn:alternate}).

If $A$ is any algebra with an action of $G$ by automorphisms, and $J<L$
are subgroups of $G$, we define the {\em transfer} map $T^L_J:A^J
\rightarrow A^L$ by $$T^L_J(a)=\sum_{h\in [L/J]} {}^h a,$$ where
$[L/J]$ is a set of representatives of the cosets $L/J$.
To prove the following theorem, we use the theory of Green functors
applied to this setting of a group action on an algebra.

\vspace{2ex}
\begin{thm} \label{prop:formula}
The cup product on $\HHD(S(V)\# G)$ induces the following
product {$\dot \smile$} on
$\bigoplus_{g\in {\mathcal C}} \HHD(S(V),S(V)\otimes g)^{Z(g)}$
under the isomorphism (\ref{eqn:alternate}):\\
For $\alpha\in \HHD(S(V),S(V)\otimes g)^{Z(g)}$ and
$\beta\in \HHD(S(V),S(V)\otimes h)^{Z(h)}$,
$$
  \alpha\ \dot\smile\ \beta 
  \ =\ \rule[-4ex]{0ex}{8ex}
\sum_{x\in D}\ T^{Z(k)}_{\rule[.5ex]{0ex}{1.5ex}
{}^yZ(g)\,\cap\, {}^{yx}Z(h)}
  ({}^y\alpha\smile {}^{yx}\beta),
$$
where $D$ is a set of representatives of the double cosets $Z(g)
\backslash G/Z(h)$, and $k=k(x)$ in $\mathcal C$ and $y=y(x)$ are chosen so that
$\ {}^y g\cdot\, {}^{yx}h = k$.
The product $ {}^y\alpha\smile {}^{yx}\beta$ is the cup product 
in  $\HHD(S(V),S(V)\# G)^{ {}^yZ(g)\cap {}^{yx}Z(h)}$,
to which
we apply the transfer map $T^{Z(k)}_{\rule[.5ex]{0ex}{1ex}
{}^yZ(g)\,\cap\, {}^{yx}Z(h)}$
to obtain an element in the $k$-component $\HHD(S(V), S(V)\ot k)^{Z(k)}$.
\end{thm}
\vspace{1ex}
\begin{proof}
Let $A=\HHD(S(V),S(V)\# G)$ (considered as a $G$-algebra).
We obtain a standard Green functor by sending a subgroup $L$ of $G$
to the invariant subring $A^L$. The restriction maps (of the 
functor) are the inclusions, and the transfer maps are as defined above.
The component $\HHD(S(V),S(V)\otimes g)^{Z(g)}$ is contained in
$A^{Z(g)}$. In the decomposition (\ref{eqn:alternate}), 
$\alpha\in \HHD(S(V),S(V)\otimes g)^{Z(g)}$ on the right side is
identified with $T^G_{Z(g)}(\alpha)$ on the left side, and similarly
for $\beta$.
The formula in the proposition is the standard one for the
product of $T^G_{Z(g)}(\alpha)$ and $T^G_{Z(h)}(\beta)$ given by
the Mackey formula
(e.g., see~\cite[Prop.\ 1.10]{Thevenaz}, due to Green). 
\end{proof}
\vspace{2ex}

The product formula in
Theorem~\ref{prop:formula} is more than a theoretical observation,
but useful in computations, as the next example shows:

\vspace{3ex}
\begin{example}\label{C6S3}
Let $V=\CC^6$ and $G=\Sym_3$.
Let $v_1,w_1,v_2,w_2,v_3,w_3$ be the standard orthonormal basis of $\CC^6$,
and let $G$ act on $V$ via the permutation representation on $\{v_1,v_2,v_3\}$
and on $\{w_1,w_2,w_3\}$.
Choose $1, (12), (123)$ as representatives of the conjugacy classes
of $\Sym_3$. By  (\ref{eqn:alternate}),
$
  \HHD(S(V)\#\Sym_3)
$
is isomorphic to $$ \HHD(S(V))^{\Sym_3}\oplus
  \HHD(S(V),S(V)\otimes (12))^{\langle (12)\rangle} \oplus 
   \HHD(S(V),S(V)\otimes (123))^{\langle (123)\rangle}.
$$
These summands are as follows, considering
that the actions of $(12)$ and of $(123)$ on the latter two summands,
respectively, are trivial:
{\begin{small}
\begin{eqnarray*}
\HHD(S(V))^{\Sym_3} \! & = & \! (S(V)\ot \Wedge^{\DOT}V^*)^{\Sym_3}\\
\HHD(S(V),S(V)\otimes (12))^{\langle (12)\rangle}\! &=& \!
   S(V^{(12)})\ot \Wedge^{\DOT -2} (V^{(12)})^*\ot
  (v_1^*-v_2^*)\wedge (w_1^*-w_2^*)\ot (12)\\
\HHD(S(V),S(V)\otimes (123))^{\langle (123)\rangle} \! &=& \! 
  S(V^{(123)})\ot \Wedge^{\DOT -4} (V^{(123)})^*\ot \\
 && 
  (v_1^*-v_2^*)\wedge (w_1^*-w_2^*)\wedge
  (v_2^*-v_3^*)\wedge (w_2^*-w_3^*)\ot (123)  
\end{eqnarray*}
\end{small}}\noindent
where the fixed point spaces are
$V^{(12)} = \Span_{\CC} \{v_1+v_2 , \ w_1+w_2, \ v_3, \ w_3 \}$
and $V^{(123)} = \Span_{\CC}\{v_1+v_2+v_3, \ w_1+w_2+w_3\}$.
The product of an element of $\HHD(S(V))^{\Sym_3}$ with an element in
any of the three components is given simply by the
componentwise product on the exterior algebras and the symmetric
algebras. The product of any element of the $(123)$-component with
any element of the $(12)$-component is 0 as the exterior product is
of linearly dependent elements. 
The product of any two elements of the $(123)$-component
is 0 by degree considerations. 

It remains to determine products
of pairs of elements from the $(12)$-component. 
In the notation of Theorem~\ref{prop:formula}, $g=(12)=h$, and
we may take the set $D$ of representatives of the double cosets
$\langle (12)\rangle \backslash \Sym_3 /\langle (12)\rangle$ to
be $\{1, (123)\}$. If $x=1$, we have $(12)\cdot (12) = 1$, 
so we take $y=1$; 
in any corresponding  product $\alpha\smile\beta$, 
the exterior part is a product of linearly dependent
elements. Thus the product corresponding to  this choice of $x$
must be 0. Now consider the case when $x=(123)$. We have
$(12)\cdot {}^{(123)}(12)=(12)(123)(12)(132) =(123)$, and so we may
take $y=1$. 

We now have, for any pair $\alpha, \beta\in \HHD(S(V),S(V)\ot(12))$, 
the corresponding cup product in $\HHD(S(V)\# \Sym_3)$
given in terms of the decomposition~(\ref{eqn:alternate}):
$$ 
  \alpha\, \dot\smile\, \beta=\  T_{\{1\}}^{\langle\,
  (123)\,\rangle} (\alpha\smile {}^{(123)}\beta ).
$$
For example, 
$$\vol^{\perp}_{(12)}
\,\dot\smile\, \vol^{\perp}_{(12)} =  T_{\{1\}}^{\langle\, (123)\,\rangle }
(\vol^{\perp}_{(12)}\, \smile\, \vol^{\perp}_{(23)}) = 3 \vol^{\perp}_{(123)}.$$
The efficient formula of Theorem~\ref{prop:formula} saves time
by reducing nine
product calculations to one in this case.
\end{example}
\vspace{3ex}

\vspace{5ex}
{\sc Acknowledgements.}
The second author thanks Ragnar-Olaf Buchweitz for many helpful conversations,
particularly regarding an early version of Theorem~\ref{cupformula}.
The authors thank Victor Reiner for aid with Lemma~\ref{G(r,1,n)} and for
finding Example~\ref{Victor} for them.



\end{document}